\documentclass[10pt,leqno]{article}

\usepackage{a4}
\usepackage{times}

\def\int{\mathbb{Z}}
\def\e{\varepsilon}
\def\Ue{{\cal U}_{\varepsilon}({\mathfrak g})}
\def\g{\mathfrak g}
\def\Spec{\text{Spec}~}
\def\O{{\cal O}}
\def\G{{\bar G}}
\def\l{{\mathfrak l}}\def\h{{\mathfrak h}}
\def\k{{\mathfrak k}}
\def\p{{\mathfrak p}}

\def\gg{{\bar{\g}}}
\def\proof{{\bf Proof. }}
\def\s{\dot{s}}
\def\w{\dot{w}}
\def\z{\dot{z}}
\def\wp{well-placed}
\def\Pf{\proof}
\def\ph{\p^u/\h^u}
\def\K{K^\circ}

\usepackage{graphics}
\usepackage{amssymb}
\usepackage{amsmath}
\input xy
\xyoption{all}
\title{Spherical orbits and representations of $\Ue$}
\author{{\sc N.\ Cantarini}\thanks{{\tt cantarin$@$math.unipd.it}}~~~
{\sc G.\ Carnovale}\thanks{{\tt carnoval$@$math.unipd.it}}~~~
{\sc M.\ Costantini}\thanks{{\tt costantini$@$math.unipd.it}}
}

\newtheorem{theorem}{Theorem}[section]
\newtheorem{lemma}[theorem]{Lemma}
\newtheorem{corollary}[theorem]{Corollary}
\newtheorem{proposition}[theorem]{Proposition}
\newtheorem{definition}[theorem]{Definition}
\newtheorem{remark}[theorem]{Remark}

\date{}
\begin{document}
\maketitle
\date{}
\begin{abstract}
Let $\Ue$ be the simply connected quantized enveloping algebra at roots
of one
associated to a finite dimensional complex simple Lie algebra $\mathfrak g$.
The De Concini-Kac-Procesi conjecture on the dimension of the
irreducible representations of $\Ue$  is
proved for the representations corresponding to the spherical conjugacy
classes of the simply connected algebraic group $G$ with Lie algebra $\mathfrak g$. We achieve this result by means of a new characterization of the spherical conjugacy classes of $G$ in terms of elements of the Weyl group.

\end{abstract}
\section*{Introduction}
Since their appearance in the mid 80's quantum groups have been
extensively investigated. In particular the representation theory of the
quantized enveloping algebra ${\cal U}_\varepsilon(\mathfrak g)$, as introduced
in \cite{DC-K}, and of the quantum function algebra $F_\varepsilon[G]$
(\cite{DC-L}) has been deeply studied by many authors. Here $\mathfrak g$
is a simple complex Lie algebra, $G$ is the
corresponding simple simply-connected algebraic group, and $\varepsilon$
is a primitive $\ell$-th root of unity, with $\ell$ an odd integer strictly
greater than $1$. However, while the irreducible
representations of $F_\varepsilon[G]$ are well described (\cite{DC-P}),
the representation theory of ${\cal U}_\varepsilon(\mathfrak g)$ is far from
being understood. In this context there is a procedure to associate
a certain conjugacy class $\O_V$   of $G$ to each simple ${\cal
U}_\varepsilon(\mathfrak g)$-module $V$. The De Concini-Kac-Procesi
conjecture asserts that $\ell^{\frac{1}{2} \dim \O_V}$ divides $\dim V$. At
present the conjecture has been proved only in some cases, namely for
the conjugacy classes of maximal dimension - the regular orbits -
(\cite{DC-K-P2}), for the subregular
unipotent orbits in type $A_n$ when $\ell$ is a power of a prime
(\cite{Nico3}),  for all orbits in $A_n$ when $\ell$ is a
prime (\cite{Nico1}), and for the conjugacy classes $\O_g$ of $g\in SL_n$
when the conjugacy class of the unipotent part of $g$ is spherical
(\cite{Nico2}). We recall that a conjugacy class $\O$ in $G$ is
called {\it spherical} if there exists a Borel subgroup of $G$ with a
dense orbit in $\O$.
The proof of the conjecture in \cite{Nico2} makes use of the
representation theory of the quantized Borel
subalgebra $B_{\e}$ introduced in \cite{3}. This method
works for representations corresponding to unipotent spherical
orbits and this underlines the correspondence between the geometry of
the conjugacy class and the structure of the corresponding irreducible
representations.

The same approach is extended in the present paper to the case of any
simple Lie algebra $\mathfrak g$ and any spherical
conjugacy class of $G$.
For this purpose we make use of the analysis of the spherical conjugacy
classes in $G$. In order to
determine the semisimple ones we use the classification of spherical pairs
$(G,H)$ where $H$ is a
closed connected reductive subgroup of $G$ of the same rank (see \cite{kramer},
\cite{bri}). On the other hand
the spherical unipotent conjugacy classes (or equivalently the  spherical
nilpotent adjoint orbits in
$\mathfrak g$) have been classified by  Panyushev in \cite{pany}
(see also \cite{pany4} 
for a proof which does not
rely on the classification of nilpotent orbits). We finally determine the
remaining spherical
conjugacy classes in Section 2.3.

Our strategy in the proof of the De Concini-Kac-Procesi conjecture for
representations corresponding to spherical orbits relies on a so far
unknown  characterization of these orbits in terms of elements of the Weyl
group $W$ of $G$. More
precisely, let us fix a pair of opposite Borel subgroups $(B,B^-)$. If $\O$
is any conjugacy class in
$G$, there exists a unique element $z=z(\O)\in W$ such that $\O\cap
B\z B$ is
open dense in $\O$.
We give a  characterization of spherical conjugacy classes  in the
following theorem:

\medskip

{\bf Theorem 1.}{\it ~Let $\O$ be a conjugacy class in   $G$, $z=z(\O)$.  Then
$\O$ is spherical
if and only if
$\dim \O = \ell(z)+rk(1-z)$.}

\medskip

Here $\ell(z)$ denotes the length of $z$ and $rk(1-z)$ denotes the rank of $1-z$ in the standard representation of
$W$. In order to
make use of
the  representation theory of  $B_{\e}$, we show that if $\O$ is a
spherical conjugacy class, then
$\O\cap B\z(\O)B\cap B^-$ is always non-empty. As a consequence of this fact
we obtain our main result
on the representation theory of $\Ue$:

\medskip

{\bf Theorem 2.}{\it ~Assume  $\mathfrak{g}$ is a finite dimensional
simple
complex Lie algebra and $\ell$ is a good integer. If $V$ is a simple ${\cal
U}_\varepsilon(\mathfrak{g})$-module whose associated conjugacy class
$\O_V$ is spherical, then
$\ell^{\frac{1}{2} \dim \O_V}$ divides $\dim V$.}

\medskip

The paper is structured as follows. In Section 1 we introduce notation
and recall the
classification of the spherical nilpotent orbits 
of $\g$.
In Section 2 the spherical conjugacy classes of $G$ are analyzed and
the main theorems are proved. In establishing Theorem 1 we shall
deal with the classical and the exceptional cases separately and we
shall consider first the unipotent conjugacy classes of $G$, then the
semisimple conjugacy classes and, finally, the conjugacy classes of
$G$ which are neither unipotent nor semisimple. 
Section \ref{conclusions} is dedicated to the analysis of the
properties of the correspondence $\O\longmapsto z(\O)$ when $\O$ is a
spherical conjugacy class.   In
Section 3 the De Concini-Kac-Procesi conjecture is proved for
representations corresponding to spherical conjugacy classes. The
proof is then
extended, using the De Concini-Kac reduction theorem (\cite{reduction}), to a larger
class of representations (see Corollary \ref{larger}). As a consequence,
the De
Concini-Kac-Procesi conjecture is proved in type $C_2$.

As far as notation and terminology are concerned we shall follow
\cite{DC-K} and \cite{Hu}. In particular for the definition of the
classical groups we choose the bilinear forms associated to the
following matrices with respect to the canonical bases:
$$\left(\begin{array}{cc}
0 & I_n\\
-I_n & 0
\end{array}\right) ~\mbox{for}~ C_n,~
\left(\begin{array}{cc}
0 & I_n\\
I_n & 0
\end{array}\right) ~\mbox{for}~ D_n,~
\left(\begin{array}{ccc}
1 & 0 & 0\\
0 & 0 & I_n\\
 0 & I_n & 0
\end{array}\right) ~\mbox{for}~ B_n.$$
In each case we fix the Borel subgroup corresponding to the set of
simple roots as described in \cite[\S 12.1]{Hu}.

{\center\bf Acknowledgements.}
The authors would like to thank Andrea Maffei for helpful
discussions and suggestions. 

\section{Preliminaries} 
Let us introduce the objects of our investigation. 
Let $A$ be an $n\times n$ Cartan matrix and let ${\mathfrak g}$ be the
associated simple complex Lie algebra, with Cartan subalgebra
${\mathfrak h}$. Let $\Phi$ be the  set of roots relative to
${\mathfrak h}$, $\Phi^+$ a fixed set of positive roots and
$\Delta=\{\alpha_1, \dots, \alpha_n\}$
the corresponding set of simple roots. Let $G$ be a reductive
algebraic group with Lie algebra ${\mathfrak g}$,  $T$  the maximal
torus with Lie algebra ${\mathfrak h}$, $B$ the Borel
subgroup
determined by $\Phi^+$ and $B^-$ the Borel subgroup opposite to $B$. Let $U$
(resp.\ $U^-$) be the unipotent radical of $B$ (resp.\ $B^-$). 

Let $W$ be the Weyl group of ${\mathfrak g}$ and let us denote by
$s_{\alpha}$ the reflection corresponding to the root $\alpha$. By
$\ell(w)$ we shall denote the length of the element $w\in W$ and
by $rk(1-w)$ we shall mean the rank of $1-w$ in the
standard representation of the Weyl group. By $w_0$ we shall denote
the  longest element in $W$. 
If $N=N(T)$ is the normalizer
of $T$ in $G$  then $W=N/T$; given an element $w\in W$ we shall denote
a representative of $w$ in $N$ by $\dot{w}$.
For any root $\alpha$ of $\g$ we shall denote by $x_{\alpha}(t)$ the
elements of the corresponding root subgroup $X_{\alpha}$ of $G$. We
shall choose the representatives $\dot{s}_{\alpha}\in N$ of the reflection
$s_{\alpha}\in W$ as in \cite[Theorem 7.2.2]{Carter1}.
In particular we recall that the Weyl group of $Sp_{2n}$ (resp.\ $SO_{2n}$) can be identified
with the group of permutations $\sigma$ in the symmetric group
$S_{2n}$ such that $\sigma(n+i)=\sigma(i)\pm n$ for all $1\leq i\leq
n$ (resp.\ $\sigma(n+i)=\sigma(i)\pm n$ and $\#\{i\leq n~|~\sigma(i)>n\}$ is
even) and it is exactly for these elements that one can choose a monomial
representative in  $Sp_{2n}$ (resp.\ $SO_{2n}$). For further details see
 \cite[p.\ 397]{F-H}.
In case of ambiguity  
we will denote  the  Weyl group (resp.\  Borel subgroups) of an
algebraic group $K$ by $W(K)$ (resp.\ $B(K)$, $B^{-}(K)$).

In order to describe the unipotent  conjugacy classes of $G$ we will make use
of their standard descriptions in terms of Young diagrams and weighted
Dynkin diagrams  \cite[\S 13.1, \S 5.6]{2}. 
For the dimension of these classes  we will
refer to \cite[\S 13.1]{2}.

%

\begin{definition}\label{sferica} Let $K$ be a connected algebraic group over $\mathbb C$ and let $H$ be a closed subgroup of $K$.
The homogeneous space $K/H$ is called spherical if  there exists a Borel subgroup of $K$ with  a dense orbit.
\end{definition}
Let us recall that the sphericity of $K/H$ 
depends only on the Lie algebras
of $K$ and $H$. By an abuse of
notation, in order to lighten the presentation, we shall identify
isogenous groups whenever convenient.

\medskip

If $\g$ is of classical type its spherical nilpotent orbits  are classified in
the following theorem: 
\begin{theorem}\label{pany}\cite[\S 4]{pany} 
The spherical nilpotent orbits in type $A_n$ and $C_n$ are those corresponding to Young 
diagrams with at most two columns. The spherical nilpotent
orbits in type $B_n$ and $D_n$ are those corresponding to Young
diagrams with at most two columns or to Young diagrams with three
columns and only one row with three boxes.\end{theorem}

In order to deal with the exceptional Lie algebras we shall also make
use of the following theorem:
\begin{theorem}\label{balatipo}\cite[Theorem 3.2]{pany2} The spherical nilpotent orbits in $\g$ are those of  type $rA_1+s\tilde{A}_1$. 
\end{theorem}

\section{Spherical conjugacy classes}
\begin{definition}\label{lyingover}
We say that an element $x\in G$ lies over an element $w\in W$ if
$x\in B\w B$.
\end{definition}

Let $\O$ be a conjugacy class in $G$. There exists a unique
element $z=z(\O)
\in W$
such that $\O\cap B\dot{z}B$ is open dense in $\O$. In particular 
\begin{equation}
\overline{\O}=\overline {\O\cap
B\z B}\subseteq \overline {B\z B}.
\label{inclusion}
\end{equation}

It follows that if $y$ is an element of $\overline\O$ and $y\in B\w B$, then
$w\leq z$ in the
Chevalley-Bruhat order of $W$. 

Let us observe that if $\O$ is a spherical conjugacy class of $G$ and
if $B.x$ is the dense $B$-orbit in $\O$, then $B.x\subseteq B\dot{z}B$.
 
\begin{theorem}\label{metodo}
Suppose that $\O$ contains an element $x\in B\dot{w}B$. 
Then $$\dim B.x
\geq \ell(w)+rk(1-w).$$ In
particular $\dim {\O} \geq \ell(w)+rk(1-w)$.
If, in addition, $\dim\O\leq\ell(w)+rk (1-w)$ then $\O$ is spherical, 
$w=z(\O)$ and $B.x$ is the dense $B$-orbit in $\O$.
\end{theorem}
\proof Let $U^w=U\cap \dot{w}U^-\dot{w}^{-1}$ and let
$B^w=U^wT$. Let us estimate the dimension of the orbit $B^w. x$. 

\medskip

{\em Step 1.} The centralizer $C_{B^w}(x)$ is contained in a maximal torus.

\noindent
Let $x=\bar{u}\dot{w}b$ be the unique decomposition of $x$ in
$U^w\dot{w}B$ and 
let $u$ be a unipotent element in $C_{B^w}(x)$. Then
$$u\bar{u}\dot{w}b=ux=xu=\bar{u}\dot{w}bu.$$
By the uniqueness of the decomposition it follows that $u=1$,
since $bu\in B$ and $u\in U^w$. Therefore the unipotent radical of
$C_{B^w}(x)$ is trivial and, by \cite[Proposition 19.4(a)]{Hu2},
$C_{B^w}(x)$ is contained in a maximal torus.

\medskip

{\em Step 2.} We have: $\dim C_{B^w}(x)\leq n-rk(1-w)$.

\noindent
Without loss of generality we may assume that $C_{B^w}(x)$ is contained in $T$. Let
$t\in C_{B^w}(x)$.  Then
$xtx^{-1}=t$ and, by \cite[\S 3.1]{spst},
$\dot{w}t\dot{w}^{-1}=t$. Therefore $C_{B^w}(x)\subset T^w$, where
$$T^w=\{t\in T ~|~ \dot{w}t\dot{w}^{-1}=t\},$$
thus
$\dim C_{B^w}(x)\leq \dim
T^w=n-rk(1-w)$. 

\medskip
Now let us observe that:
$$\dim B^w.x=\dim B^w-\dim
C_{B^w}(x)=$$
$$=\ell(w)+n-\dim C_{B^w}(x)\geq
\ell(w)+n-n+rk(1-w)=\ell(w)+rk(1-w).$$
It follows that if, in addition, $\ell(w)+rk(1-w)\geq \dim\O$,
then $\dim C_{B^w}(x)=\dim T^w$ and $\dim\O=\ell(w)+rk(1-w)$.
In particular $B.x$ is the dense $B$-orbit in $\O$.\hfill$\Box$
         
\begin{proposition}\label{metodo2} Let $\O$ be a conjugacy class, $z=z(\O)$. 
If there exists an
element $w\in W$ such that $w\leq z$ and $\dim \O \leq \ell(w)+rk(1-w)$,
then $\O$ is spherical with $\dim\O=\ell(w)+rk(1-w)=\ell(z)+rk(1-z)$.
\end{proposition} 
\Pf From $w\leq z$ it follows that $\ell(w)+rk(1-w) \leq
\ell(z)+rk(1-z)$. Indeed it is
enough to consider the case $\ell(z)=\ell(w)+1$: then $rk(1-z)=rk(1-w)\pm
1$ so that either
$\ell(z)+rk(1-z)=\ell(w)+rk(1-w)+2$ or $\ell(z)+rk(1-z)=\ell(w)+rk(1-w)$
and the inequality follows.
Therefore $\dim \O\leq \ell(w)+rk(1-w) \leq \ell(z)+rk(1-z)$.
By Theorem \ref{metodo}  we obtain
$$
\dim \O=\ell(z)+rk(1-z) = \ell(w)+rk(1-w).
$$
\hfill$\Box$

Let us observe that  it may happen that $w\not=z$. 

\begin{corollary}\label{1} Let $\O$ be a conjugacy class, $z=z(\O)$. 
Let $w_1,\ldots,w_k$
be elements of $W$ such that $\O\cap B\dot{w}_iB\not=\emptyset$ for
$i=1,\ldots,k$, and let us consider the
set $X$ of minimal elements in
$$
\{w\in W\mid w\geq w_i, ~i=1,\dots,k\}.
$$
If for every $w\in X$ we have $\dim \O\leq\ell(w)+rk(1-w)$, then $\O$ is spherical.
\end{corollary}
\Pf Since $w_i\leq z$ for
$i=1,\ldots,k$, there exists $w\in X$
such that $w\leq z$. Then we conclude by Proposition \ref{metodo2}. \hfill$\Box$
 
\begin{corollary}\label{reticolo} Let $\O$ be a conjugacy class. Let $w_1$, $w_2$
be elements of $W$ such that $\O\cap B\dot{w}_iB\not=\emptyset$ for $i=1$,
$2$.
If $$
\{w\in W\mid w\geq w_i, ~i=1,2\}=\{w_0\}
$$
then $z(\O)=w_0$. If, in addition, $\dim\O\le
\ell(w_0)+rk(1-w_0)$ then $\O$ is spherical. 
\hfill$\Box$
\end{corollary}

\begin{definition}\label{wp} Let $\O$ be a  conjugacy class. We say
that $\O$ is {\em well-placed} if there exists an element  $w\in W$
such that $$\O\cap B^-\cap B\dot{w}B\not=\emptyset ~~\mbox{and} ~~\dim \O=\ell(w)+rk(1-w).$$
\end{definition}

It follows from Definition \ref{wp} and Theorem \ref{metodo} that if a conjugacy class $\O$ is \wp\
then it is spherical and $z(\O)=w$.
Our aim is to show that every spherical conjugacy class is \wp. 

\medskip

In the
sequel we
will make use of
following lemma:
\begin{lemma}\label{isogeny} Let $\phi: G_1 \longrightarrow G_2$ be an isogeny of
reductive algebraic groups. Let $x_1\in G_1$,
$x_2=\phi(x_1)$ and let $\O_{x_i}$ be the conjugacy class of
$x_i$ in $G_i$. Let $w\in
W=W(G_1)=W(G_2)$ and let $\dot{w}_i$ be a representative of $w$ in
$G_i$. Then $B(G_1)\dot{w}_1B(G_1)\cap B^-(G_1)\cap
\O_{x_1}\neq\emptyset$ if and only if $B(G_2)\dot{w}_2B(G_2)\cap B^-(G_2)\cap
\O_{x_2}\neq\emptyset$.\hfill$\Box$
\end{lemma}

\subsection{Unipotent conjugacy classes}
In view of Definition \ref{wp} we begin this section with a result
concerning the intersection between $U^-$ and the
(unique) dense $B$-orbit in a spherical unipotent conjugacy class.

\begin{lemma}\label{B-unipotent} Let $\O$ be a unipotent spherical conjugacy class, $B.x$ the
(unique) dense $B$-orbit
in $\O$. Then $B.x\cap U^-$ is not empty.
\end{lemma} 
\proof Let $g\in{\O}$ and let $P$ be the canonical parabolic subgroup of $G$
associated to $g$ (see \cite{spst} and \S 2.3.1).
Then $g$ lies in
the unipotent
radical $P^u$ of $P$, and $H=C_G(g)\leq P$. Since ${\O}$ is spherical, there exists
a Borel subgroup $B_1$
of $G$ such that $HB_1$ is dense in $G$. In particular, $PB_1$ is dense
in $G$. Without loss of
generality, we may assume $P\geq B$, $P=P_{J_1}$ with $J_1\subseteq \{\alpha_1,\,\ldots,\,\alpha_n\}$ say,
and $B_1=\dot{\tau}B\dot{\tau}^{-1}$, with
$\dot{\tau}\in N_{J_1,\emptyset}$, following the notation in \cite[\S 2.8]{2} (here
we have $J_2=\emptyset$).
In our case the subset $K$ of $\{\alpha_1,\,\ldots,\,\alpha_n\}$ is empty. We recall that
$N_{J_1,\emptyset}=\{\dot \sigma\mid \sigma\in
D_{J_1,\emptyset}\}$ and that $D_{J_1,\emptyset}=D_{J_1}^{-1}$, where
$D_{J_1}=\{\sigma\in W\mid
\sigma(\Phi_{J_1}^+)\subseteq \Phi^+\}$. Then
$\tau^{-1}(\Phi_{J_1}^+)\subseteq \Phi^+$.
We show that $P\,\dot{\tau}B\dot{\tau}^{-1}$ is dense in $G$ if and only if
$\tau^{-1}(\Phi^+\setminus\Phi_{J_1})\subseteq
\Phi^-$ (which then implies that $w_0\tau$ is the longest element of
$W_{J_1}$). We have
$$
P\cap\dot{\tau}B\dot{\tau}^{-1}=(P^u\cap\dot{\tau}U\dot{\tau}^{-1})(L_{J_1}\cap
\dot{\tau}B\dot{\tau}^{-1})
$$                       
and $L_{J_1}\cap\dot{\tau}B\dot{\tau}^{-1}=L_{J_1}\cap B_1$ is a Borel subgroup of $L_{J_1}$ by
\cite[Propositions 2.8.7, 2.8.9]{2}. Let us denote by $r$ the number of positive roots in $\Phi_{J_1}$
and by $s$ the dimension of
$P^u\cap \dot{\tau}U\dot{\tau}^{-1}$. Then $P\,\dot{\tau}B\dot{\tau}^{-1}$ is dense in $G$ if and only if
$\dim(P\cap\dot{\tau}B\dot{\tau}^{-1})=\dim P+\dim
B-\dim G$. Since $\dim P=\dim P^u+\dim L_{J_1}=N+n+r$, $\dim B=N+n$, $\dim
L_{J_1}\cap B=n+r$, we
get that $P\dot{\tau}B\dot{\tau}^{-1}$ is dense in $G$ if and only $s=0$, that is $P^u\cap\dot{\tau}U\dot{\tau}^{-1}
=\{1\}$. This in
turn is equivalent to $(\Phi^+\setminus \Phi_{J_1})\cap
\tau(\Phi^+)=\emptyset$, that is
$\tau^{-1}(\Phi^+\setminus\Phi_{J_1})\subseteq
\Phi^-$, as we wanted.
 
We are now in the position to exhibit an element in $B.x\cap U^-$. By
hypothesis we have $g\in
P^u=\prod _{\beta\in\Phi^+\setminus\Phi_{J_1}}X_\beta$. Then $\dot{\tau}^{-1}g\dot{\tau}$
lies in
$\prod _{\beta\in\Phi^+\setminus\Phi_{J_1}}X_{\tau^{-1}\beta}\leq U^-$. On the
other hand, from $H\,\dot{\tau}B\dot{\tau}^{-1}$ dense in $G$ it follows that
  $C_G(\dot{\tau}^{-1}g\dot{\tau})\,
B$ is dense in $G$, hence $\dot{\tau}^{-1}g\dot{\tau}$ lies in $B.x$. \hfill$\Box$

\medskip

Let us observe that we can directly deal with the minimal unipotent
conjugacy class:
\begin{proposition}\label{minima} Let $\O$ be the unipotent conjugacy class of type $A_1$ (minimal orbit). Then $\O$ is \wp.
\end{proposition} 
\Pf Let $\beta_1$ denote the highest root of $\g$. Then
$x_{-\beta_1}(1)$  is a representative  of $\O$. For every
positive root $\alpha$ and every $t\neq 0$ we have:
\begin{equation}
x_{-\alpha}(t)=x_{\alpha}(t^{-1})h\dot{s}_{\alpha}x_{\alpha}(t^{-1})
\label{reflection}
\end{equation}
for some $h\in T$ (see \cite[p.\ 106]{Carter1}).
In particular $x_{-\beta_1}(1)$ belongs to $B\dot{s}_{\beta_1}B\cap B^-$. By
\cite[Lemma 4.3.5]{CMG} we have 
$$\ell({s}_{\beta_1})+rk(1-{s}_{\beta_1})=\#\{\alpha\in\Phi^+\ |\
\alpha\not\perp\beta_1\}+1=\dim\O$$ and the statement follows. \hfill$\Box$

\subsubsection{Classical type}
This section is devoted to the analysis of the spherical unipotent
 conjugacy classes of $G$ when  $G$ is of classical type. Since the
 case of type $A_n$ has been treated in \cite{Nico2} we shall assume
 that $G$ is of type $B_n$, $C_n$ or $D_n$.

It will be useful
for our purposes to fix some notation for  Young diagrams corresponding to 
spherical unipotent conjugacy classes. We will denote by 
$X_{t,m}$ and $Z_{t,m}$, respectively,   the following Young diagrams with $m$ boxes:
\begin{equation}\label{diagrammi}
\begin{picture}(30,105)


\put (-60,100){\line(1,0){20}}
\put (-60,90){\line(1,0){20}}
\put (-60,70){\line(1,0){20}}
\put (-60,60){\line(1,0){20}}
\put (-60,50){\line(1,0){10}}
\put (-60,40){\line(1,0){10}}
\put (-60,10){\line(1,0){10}}
\put (-60,0){\line(1,0){10}}


\put (-60,100){\line(0,-1){10}}
\put (-60,70){\line(0,-1){30}}
\put (-60,10){\line(0,-1){10}}
\put (-50,100){\line(0,-1){10}}
\put (-50,70){\line(0,-1){30}}
\put (-50,10){\line(0,-1){10}}
\put (-40,100){\line(0,-1){10}}
\put (-40,70){\line(0,-1){10}}

\put (-100, 50){$X_{t,m}=$}
\put (10, 50){$Z_{t,m}=$}
\put (-51,75){\vdots}
\put (-56,20){\vdots}
\put (-36,76){$\Biggr\}$}
\put (-26,77){$t$}





\put (50,100){\line(1,0){30}}
\put (50,90){\line(1,0){30}}
\put (50,80){\line(1,0){20}}
\put (50,60){\line(1,0){20}}
\put (50,50){\line(1,0){20}}
\put (50,40){\line(1,0){10}}
\put (50,30){\line(1,0){10}}
\put (50,10){\line(1,0){10}}
\put (50,0){\line(1,0){10}}


\put (50,100){\line(0,-1){10}}
\put (60,100){\line(0,-1){10}}
\put (70,100){\line(0,-1){10}}
\put (80,100){\line(0,-1){10}}
\put (50,90){\line(0,-1){10}}
\put (50,60){\line(0,-1){30}}
\put (50,10){\line(0,-1){10}}
\put (60,90){\line(0,-1){10}}
\put (60,60){\line(0,-1){30}}
\put (60,10){\line(0,-1){10}}
\put (70,90){\line(0,-1){10}}
\put (70,60){\line(0,-1){10}}

\put (59,65){\vdots}
\put (54,15){\vdots}
\put (64,0){.}
\put (74,66){$\Biggr\}$}
\put (84,67){$t$}

\end{picture}
\end{equation}
By an abuse of notation, given a unipotent element $u\in G$ and a 
Young diagram of fixed shape $J$, we will say that $u=J$ if the conjugacy
 class of $u$ is described by $J$.
 
It will be convenient for our purposes to understand when an element
of a classical group  lies over the longest element $w_0$ of the Weyl
group.
\begin{remark}\label{bigcell}{\rm
Let $G=Sp_{2n}$ (resp.\ $SO_{2n}$ and $n$ even) so that $w_0=-1$.
Then the elements 
of $B^-$ and $B$ are  of the form $\left(
\begin{array}{c|c}
^tF^{-1} & 0\\
\hline
F\Sigma & F
\end{array}
\right)$ and   $\left(
\begin{array}{c|c}
X & XA\\
\hline
0 & ^tX^{-1}
\end{array}
\right)$, respectively, where $\Sigma$ and $A$ are symmetric 
(resp.\ skew-symmetric), and $F$ and $X$ are  upper
triangular, invertible matrices. Therefore an element $x\in B^-$ lies
over $w_0$ if there exist  upper triangular invertible matrices $X$ and $Y$,
and symmetric (resp.\ skew-symmetric)
matrices $A$, $B$ such that
\begin{equation}
x=\left(
\begin{array}{c|c}
^tF^{-1} & 0\\
\hline
F\Sigma & F
\end{array}
\right)=
\left(
\begin{array}{c|c}
X & XA\\
\hline
0 & ^tX^{-1}
\end{array}
\right)\left(
\begin{array}{c|c}
0 & I_{n}\\
\hline
\mp I_{n} & 0
\end{array}
\right)\left(
\begin{array}{c|c}
Y & YB\\
\hline
0 & ^tY^{-1}
\end{array}
\right).
\label{M1}
\end{equation}

A direct computation shows that  (\ref{M1})  
holds if and only if $F\Sigma=~^t\!X^{-1}Y$, i.e., if and only if $F\Sigma$
lies in the big cell of 
$GL_n$ or, equivalently, if its principal minors  are different from
zero (see, for example, \cite[Exercise 28.8]{Hu2}).

Similarly, if  $G=SO_{2n+1}$, so that $w_0=-1$,
 the elements of $B^-$ and 
$B$ are of the form
$\left(
\begin{array}{c|c|c}
1 & ^t\psi & 0 \\
\hline
0 & ^tF^{-1} & 0  \\
\hline
-F\psi& F\Sigma & F
\end{array}
\right)$ and  
 $\left(
\begin{array}{c|c|c}
1 &0 & ^t\gamma \\
\hline
-X\gamma & X & XA\\
\hline
0& 0 & ^tX^{-1}
\end{array}
\right)$, respectively, where the symmetric parts of  $\Sigma$ and $A$ 
are $-(1/2)\psi\,^t\!\psi$ and $-(1/2)\gamma\,^t\!\gamma$
 respectively, and $F$ and $X$ are  upper
triangular, invertible matrices. Therefore an element $x=\left(
\begin{array}{c|c|c}
1 & ^t\psi & 0 \\
\hline
0 & ^tF^{-1} & 0  \\
\hline
-F\psi& F\Sigma & F
\end{array}
\right)\in B^-$ lies over 
$w_0$ if and only if there exist two upper triangular invertible matrices 
$U$ and $X$, two vectors $\gamma$ and $c$, and two matrices $A$ and $S$ 
with symmetric part equal to $-(1/2)\gamma\,^t\!\gamma$ and
 $-(1/2)c\,^t\! c$, respectively, such that the following equality holds:
\begin{equation}\label{M2}
x=
\left(
\begin{array}{c|c|c}
1 &0 & ^t\gamma \\
\hline
-X\gamma & X & XA\\
\hline
0& 0 & ^tX^{-1}
\end{array}
\right)
\left(
\begin{array}{c|c|c}
(-1)^n & 0 & 0\\
\hline
0& 0 & I_{n}\\
\hline
0& I_{n} & 0
\end{array}
\right)
\left(
\begin{array}{c|c|c}
1 &0 & ^tc \\
\hline
-Uc & U & US\\
\hline
0& 0 & ^tU^{-1}
\end{array}
\right).
\end{equation}

A tedious but straightforward computation shows that 
(\ref{M2}) holds if and only if $F\Sigma$ lies in the big cell of $GL_n$ 
and $^t\!\psi \Sigma^{-1}\psi=(-1)^n-1$. 
}\end{remark}
\begin{theorem}\label{classical-unipotent}
Let ${\cal O}_g$ be a spherical unipotent conjugacy class of an element
$g\in G$.
Then  ${\cal O}_g$ is \wp.
\end{theorem} 
\proof {\em G of type $C_n$.} 
For every integer $k=1,\dots, n$ let us consider the unipotent conjugacy class
${\cal O}_k$ of $Sp_{2n}$ parametrized by a Young diagram of
shape $X_{k,2n}$. We have: $\dim \O_k=k(2n-k+1)$. 

For every fixed $k$ let us  choose the following matrix $A_k$ in ${\cal O}_k\cap B^-$: 
$$A_k=\left (\begin{array}{l|r}
I_n & 0_n\\
\hline
I_k^{\prime} & I_n
\end{array}
\right)$$
where  $I_k^{\prime}$ is the $n\times n$  diagonal matrix 
$I_k^{\prime}=\left(\begin{array}{c|c}
I_k & 0\\
\hline
0 & 0_{n-k}
\end{array}
\right)$.

In $W$ let us consider the element
$w_k$  sending $e_i$ to
$-e_i$ for every $i=1, \dots, k$ and fixing all the other elements of
the canonical 
basis of $\mathbb{C}^n$. 
We have: $$rk(1-w_k)+\ell(w_k)=k(2n-k+1)=\dim {\cal O}_k.$$ 
If we choose the representative 
$$\dot{w}_k=\left (\begin{array}{cc|cc}
0 & 0 & I_k & 0\\
0 & I_{n-k} & 0 & 0\\
\hline
-I_k & 0 & 0 & 0\\
0 & 0 & 0 & I_{n-k}
\end{array}
\right),$$
then
$$A_k=U_k\dot{w}_kB_k$$
where $U_k=\left(\begin{array}{c|c}
I_n~~~ & I_k^{\prime}
\\
\hline
0_n~~ & I_n~~~
\end{array}
\right)$ and $B_k=\left(\begin{array}{c|c}
{\begin{array}{cc}
-I_k & 0\\
0 & I_{n-k}
\end{array}} & -I_k^{\prime}
\\
\hline
0_n & {\begin{array}{cc}
-I_k & 0\\
0 & I_{n-k}
\end{array}}
\end{array}
\right)$.

\bigskip
\noindent
This identity shows that  $A_k$ lies over $w_k$ since $U_k$ and $B_k$
belong to $B$.
This concludes the proof for $G$ of type $C_n$.

\bigskip
\noindent
{\em G of type $D_n$.} Let us consider the
 unipotent  conjugacy
classes  of $SO_{2n}$  associated to 
Young diagrams  either of shape $X_{2k,2n}$ with $k=1,\ldots,\,[n/2]$,
 or of shape $Z_{2k, 2n}$
with $k=0,\,\ldots,\,[n/2]-1$. Let us recall that when $n$ is even there are two
 distinct conjugacy classes, 
${\cal O}_{n/2}$ and ${\cal O'}_{n/2}$,
 associated to the Young diagram of shape
 $X_{n,2n}$, with 
 weighted Dynkin diagrams

$\begin{array}{cl}
{\xymatrix@R=6pt@C=18pt{
&&&&& *+[o][F]{}\ar@{-}[ld]_>>>>>>>(-.5){0}\\
D=& *+[o][F]{}\ar@{-}[r]^<<(0){0}&  *+[o][F]{}\ar@{-}[r]^<<(0){0} & \cdots  &
 *+[o][F]{} \ar@{-}[l]_<<(0){0} \\
&&&&& *+[o][F]{}\ar@{-}[lu]_>>>(.15){2}
}} & {\begin{array}{c} \\ \\~~~{\mbox{and}}\end{array}} \\
& \\
{\xymatrix@R=6pt@C=18pt{
&&&&& *+[o][F]{}\ar@{-}[ld]_>>>>>>>(-.5){2}\\
D^{\prime}=&*+[o][F]{}\ar@{-}[r]^<<(0){0} & *+[o][F]{}\ar@{-}[r]^<<(0){0} & \cdots  &
 *+[o][F]{}\ar@{-}[l]_<<(0){0} \\
&&&&& *+[o][F]{}\ar@{-}[lu]_>>>(.15){0}
}} & {\begin{array}{c} \\ \\~~~{\mbox{respectively.}}\end{array}}
\end{array}$

\medskip

Moreover let
${\cal O}_k$ (for $k<\frac{n}{2}$)  and $\tilde{\O}_k$ be the unipotent conjugacy classes
with
 Young diagrams $X_{2k,2n}$ and $Z_{2k,2n}$, respectively.

We have: $\dim {\cal
O}_k=2k(2n-2k-1)$, $\dim \tilde{\cal
O}_k=4(k+1)(n-k-1)$ and $\dim\O'_{n/2}=n^2-n=\dim \O_{n/2}$.

Now let us consider the following element $w_k$ in the Weyl group of $so_{2n}$:
$$w_k: \left\{
\begin{array}{lll}
e_i \longmapsto - e_{i+1} & \mbox{if}~ i ~{\mbox{is odd ~and}} &~1\leq i\leq 2k-1,\\
e_i\longmapsto - e_{i-1} & \mbox{if}~ i ~{\mbox{is even and}} &~2\leq i\leq
2k,\\
e_i\longmapsto e_i & \mbox{if}~ i>2k. &
\end{array}\right .$$
Then
$\ell(w_k)=4nk-4k^2-3k$ and $rk(1-w_k)=k$, therefore
$\ell(w_k)+rk(1-w_k)=\dim {\cal O}_k$. 

\noindent
Let us introduce the following matrices: 
$S_1=\left(
\begin{array}{cc}
0 & 1 \\
-1 & 0 \\
\end{array}
\right)$,
$S_k=\mbox{diag}(S_1,\ldots,\,S_1)$ of order $2k$,
%
$J_k=\left(
\begin{array}{c|c}
S_k & 0\\
\hline
0 & 0_{n-2k}
\end{array}
\right)$,
$u_k=\left(\begin{array}{c|c}
I_n & 0_n\\
\hline
J_k & I_n
\end{array}
\right)$ and 
$H_k=\left(\begin{array}{c|c}
I_n & -J_k\\
\hline
0_n & I_n
\end{array}
\right)$.\\

 Notice that $H_k\in B$ and  $u_k$ lies in  $\O_k\cap B^-$ for $k<n/2$. 
The weighted Dynkin diagram associated to $u_{n/2}$
shows that $u_{n/2}\in\O_{n/2}\cap B^-$.  Besides, the following identity of matrices holds:
$H_k\dot{w}_kH_k=u_k$ where $$\dot{w}_k= \left(\begin{array}{c|c}
{\begin{array}{cc}
0_{2k} & \\
 & I_{n-2k}
\end{array}} & J_k\\
\hline
J_k & {\begin{array}{cc}
0_{2k} & \\
 & I_{n-2k}
\end{array}}
\end{array}
\right).$$ This shows that $u_k$ lies over $w_k$ for $k=1,\dots, [n/2]$.

Let now  $n$ be even and $k=n/2$ and let us consider the automorphism $\hat\tau$
 of $SO_{2n}$
arising from 
the automorphism $\tau$ of the Dynkin diagram 
 interchanging $\alpha_{n-1}$ and $\alpha_n$. 
Then $u'_{n/2}=\hat\tau(u_{n/2})\in B^-$ is
 a representative of
 the conjugacy class $\O'_{n/2}$ associated to $D'$.
If we apply the map $\hat\tau$ to  the equality $u_{n/2}=H_{n/2}\dot{w}_{n/2} H_{n/2}$
we find that
$u'_{n/2}$ lies over ${w_{n/2}^\tau}=\tau w_{n/2}\tau\in W\subset Aut(\Phi)$.
As $\tau$ permutes simple roots, it is clear that $\ell(w_{n/2})=\ell(w_{n/2}^\tau)$.
 Therefore, 
$$\ell({w_{n/2}^\tau})+rk(1-w_{n/2}^\tau)=\ell(w_{n/2})+rk(1-w_{n/2})=
\dim\O_{n/2}=\dim \O'_{n/2}.$$
This concludes the proof for $G$ of type $D_n$
 and $\O$ a conjugacy class corresponding to a Young diagram of shape $X_{2k,2n}$ with $k\le[\frac{ n}{2}]$. 

\bigskip

Now we want to prove the statement for $\tilde\O_k$.
Let us first assume $n=2m$. Let $v_{m-1}=\left(\begin{array}{c|c}
^tF^{-1} & 0\\
\hline
F\Sigma & F
\end{array}
\right)$ where:\\[5pt]
$\bullet$~$F$ is
 the upper triangular $n\times n$ matrix with all diagonal elements
 equal to $1$, the first upper off-diagonal equal to $(-1,0,-1,\ldots,\,0,-1)$
 and zero elsewhere;\\[5pt]
$\bullet$~$\Sigma$ is the skew symmetric matrix whose
 first upper off-diagonal is   $(1, 0,1,0,\ldots,0,1)$,
 all further 
odd upper off-diagonals are equal to
 $(2,0,2,0,\ldots,0,2)$ and 
all even off-diagonals are equal to $(0,0,\ldots,0)$.\\[5pt]
\noindent
One can show that
$v_{m-1}\in\tilde\O_{m-1}$ and that
 $F\Sigma$ belongs to the big cell of $GL_{2m}$. Therefore, by Remark \ref{bigcell},
$v_{m-1}\in B^-$ lies over $w_0$. Observe that 
$$\ell(w_0)+rk(1-w_0)=n^2=4m^2=\dim
\tilde\O_{m-1}$$ and that equality holds also when $n=2$, i.e., when $SO_{2n}$
is not simple. Hence the statement is proved for $n$ even and $k=n/2
-1$. 
 
Let us consider the conjugacy class $\tilde\O_k$ for   $n$ not
necessarily even and
the   embedding $j_{2k+2}$ of $SO_{4k+4}$ into $SO_{2n}$:
$$j_{2k+2}: ~\left(
\begin{array}{c|c}
A & B\\
\hline
C & D
\end{array}\right)\longmapsto \left(
\begin{array}{c|c}
{\begin{array}{cc}
A & \\
 & I_{n-2k-2}
\end{array}} & {\begin{array}{cc}
B & \\
 & 0_{n-2k-2}
\end{array}}\\
\hline
{\begin{array}{cc}
C & \\
 & 0_{n-2k-2}
\end{array}} & {\begin{array}{cc}
D & \\
 & I_{n-2k-2}
\end{array}}
\end{array}\right).$$
 The embedded image of $v_{m-1}$, for $m=k+1$, belongs to $B^-$, it is a representative of $\tilde\O_k$ and
lies over $\eta_k=\left(\begin{array}{cc}
-I_{2k+2}&0\\
0&I_{n-2k-2}\\
\end{array}\right)$. 
One can check that $\ell(\eta_k)+rk(1-\eta_k)=\dim\tilde\O_k$ so the statement is proved for $G$ of type $D_n$.\\

\noindent
{\em G of type $B_n$.} 
Let us consider the unipotent conjugacy classes
${\cal C}_k$ and $\tilde{\cal C}_h$
of shape $X_{2k,2n+1}$ and 
$Z_{2h,2n+1}$, respectively, with $k=1,\dots, [n/2]~$ and $~h=0,\dots$,
$[(n-1)/2]$.  
We have: $\dim {\cal
C}_k=4nk-4k^2$ and $\dim \tilde{\cal
C}_h=2(h+1)(2n-2h-1)$.\\
Let us consider the following embedding of $SO_{2n}$ in $SO_{2n+1}$:
$$X\longmapsto \left(\begin{array}{cc}
1 & \\
 & X
\end{array}
\right).$$ Under this embedding a representative of an element $w\in W(SO_{2n})$ 
is mapped to a representative of an element in $W(SO_{2n+1})$.
Through the same embedding the Borel subgroup
 $B(SO_{2n})$ (resp.\ $B^-(SO_{2n})$) can be seen as a 
subgroup of $B(SO_{2n+1})$ (resp.\ $B^-(SO_{2n+1})$). The image
 of the representative $u_k$ of the class
 ${\cal O}_k\subset SO_{2n}$  is  a representative  of the class 
${\cal C}_k\subset SO_{2n+1}$, it belongs to $B^-(SO_{2n+1})$ and lies
 over $w_k$ where $w_k$ is the same as in the corresponding case of $SO_{2n}$. The length of $w_k$, viewed as an element of $W(SO_{2n+1})$, is
$\ell(w_k)=4nk-4k^2-k$ and $rk(1-w_k)=k$, therefore
 $\ell(w_k)+rk(1-w_k)=\dim{\cal C}_k$. Hence, we have the statement for
${\cal C}_k$.

Similarly, if $k\leq \left[n/2\right]-1$, the  image of the representative  $v_k$ of the class 
${\tilde\O}_k\subset SO_{2n}$ is a  representative of the class
 ${\cal\tilde C}_k\subset SO_{2n+1}$, it  belongs
 to $B^-(SO_{2n+1})$ and  lies over $\eta_k$ where $\eta_k$ is the
same as in the corresponding case of $SO_{2n}$.
If we view $\eta_k$ as an element of $W(SO_{2n+1})$
we obtain:
$$\ell(\eta_k)+rk(1-\eta_k)=4nk-4k^2-8k+4n-4+2k+2=\dim{\cal\tilde C}_k$$
so the statement holds for ${\cal\tilde C}_k$ with $k\le\left[n/2\right]-1$.

Let us now prove the statement for the classes corresponding to Young
 diagrams with no rows consisting of only one box, i.e., for ${\cal
 \tilde C}_{{\frac{n-1}{2}}}$ when $n$ is odd. In this case $$\dim
 {\cal \tilde C}_{{\frac{n-1}{2}}}=n^2+n=\ell(w_0)+rk(1-w_0).$$ Let us
 consider the element $v=\left(\begin{array}{c|c|c} 1&^t\psi&0\\
 \hline 0&I_n&0\\ \hline -\psi&\Sigma&I_n\\
\end{array}
\right)$
where $\psi=\,^t\!(1~0~\ldots~0)$ and $\Sigma$ is the $n\times n$ matrix with diagonal equal to
$(-1/2,0,\ldots,0)$, first upper off-diagonal equal to 
$(1,1,\ldots,1)$, first lower off-diagonal equal to 
$(-1,-1,\ldots,-1)$ and zero elsewhere. Then  $v$
 is a representative of ${\cal \tilde C}_{{\frac{n-1}{2}}}$. One can check 
that $\Sigma$ belongs to the big cell of $GL_n$ and that $^t\psi\,\Sigma^{-1}\psi=-2$. By Remark \ref{bigcell} we conclude the proof.
\hfill$\Box$

\subsubsection{Exceptional type}\label{E-U}
This section is devoted to the analysis of the spherical unipotent
conjugacy classes of $G$ when $G$ is of exceptional type. Let us
introduce some notation: we shall denote by $\beta_1$ the highest root
of $\g$ and, inductively, by $\beta_r$, for $r>1$, the highest root of
the root system orthogonal to $\beta_1, \dots, \beta_{r-1}$ when this is
irreducible. Similarly we shall denote by $\gamma_1$ the highest short
root of $\g$ and inductively, by $\gamma_r$, for $r>1$, the highest short root of
the root system orthogonal to $\gamma_1,\dots, \gamma_{r-1}$ when it is
irreducible.

\begin{theorem}\label{exceptional-unipotent}
Let ${\cal O}$ be a spherical unipotent conjugacy class.
Then $\O$ is \wp.
\end{theorem} 
\proof
The unipotent spherical conjugacy classes of
$G$ are those of type $rA_1+s\tilde{A}_1$.
We shall deal with the different types of orbits separately:

\medskip

\noindent
{\it Type ${A}_1$.} See Proposition \ref{minima}.

\medskip

\noindent
{\it Type $\tilde{A}_1$} ($\g$ of type $F_4$, $G_2$).
Let $\g$ be of type $G_2$. The element $x_{-\gamma_1}(1)$ is a representative of the class $\O$
of type $\tilde{A}_1$  and
lies over $s_{\gamma_1}$ by (\ref{reflection}). Therefore 
$z(\O)\ge s_{\gamma_1}$. Besides,
since $\overline\O$ contains the minimal conjugacy class, by (\ref{inclusion}) it
follows that $z(\O)\geq s_{\beta_1}$, hence, by Corollary \ref{reticolo},
$z(\O)=w_0$.
We now conclude using Lemma \ref{B-unipotent}
and noticing that
$\dim\O=8=\ell(w_0)+rk(1-w_0)$.

Let $\g$ be of type $F_4$.  
The element
  $x=x_{-\beta_1}(1)x_{-\beta_2}(1)$ is a
 representative  of the class of type $\tilde{A}_1$ as the calculation
 of its weighted Dynkin diagram shows. By (\ref{reflection}) $x$ belongs to
 $B\dot{s}_{\beta_1}
\dot{s}_{\beta_2}B$ and one can check that
 $\ell({s}_{\beta_1}{s}_{\beta_2})+rk(1-{s}_{\beta_1}{s}_{\beta_2})=22=\dim\O$. 
\medskip

\noindent
{\it Type $2A_1$} ($\g$ of type $E_6$, $E_7$, $E_8$).
The element $x_{-\beta_1}(1)x_{-\beta_2}(1)$  is a
 representative  of this class. By construction and
by (\ref{reflection}) $x_{-\beta_1}(1)x_{-\beta_2}(1)$ lies over
 ${s}_{\beta_1}{s}_{\beta_2}$. One can check that  $\ell({s}_{\beta_1}{s}_{\beta_2})+rk(1-{s}_{\beta_1}{s}_{\beta_2})=\dim\O$. 

\medskip

\noindent
{\it Type $3A_1$} ($\g$ of type $E_6$, $E_7$, $E_8$).
If $\g$ is of type $E_7$ there are two conjugacy classes of type
$3A_1$ that, following \cite{B-CII},  we shall denote by
$(3A_1)^{\prime}$, $(3A_1)^{\prime\prime}$. 
A representative of the class $(3A_1)^{\prime\prime}$ is
$x_{-\beta_1}(1)x_{-\beta_2}(1)x_{-\alpha_7}(1)$, as one can verify
by computing its weighted Dynkin diagram. Relation (\ref{reflection})
implies that $x_{-\beta_1}(1)x_{-\beta_2}(1)x_{-\alpha_7}(1)$ lies
over ${s}_{\beta_1}{s}_{\beta_2}{s}_{\alpha_7}$ since $\alpha_7$ is
orthogonal to $\beta_1$ and $\beta_2$. One can verify that 
$\ell({s}_{\beta_1}{s}_{\beta_2}{s}_{\alpha_7})+rk(1-{s}_{\beta_1}{s}_{\beta_2}{s}_{\alpha_7})=54=\dim\O$. 

In order to handle the remaining classes of type $3A_1$ we consider
subalgebras of type $D_4$ in $\g$ and the corresponding immersions of
algebraic groups. We realize the root systems of these subalgebras as
the sets of roots orthogonal to $\ker(1-w)$ where $w\in W$ is chosen
as follows:
\begin{itemize}
\item
$w_0=s_{\beta_1}s_{\beta_2}s_{\beta_3}s_{\alpha_4}$ if
$\g$ is of type $E_6$;  
\item
$s_{\beta_1}s_{\beta_2}s_{\alpha_2+\alpha_3+2\alpha_4+\alpha_5}s_{\alpha_3}$
if $\g$ is of type $E_7$;
\item
$s_{\beta_1}s_{\beta_2}s_{\beta_3}s_{\alpha_7}$
if $\g$ is of type $E_8$.
\end{itemize}
In Theorem \ref{classical-unipotent} we proved that if $\O^{\prime}$
is  the class of type $3A_1$ of
$D_4$ then $z(\O^{\prime})$ is  the longest element of the Weyl group
of $D_4$. By
construction, in each case $w$ is the longest element of the Weyl
group of the corresponding copy of
$D_4$. One can verify that  $\ell(w)+rk(1-w)$ is
equal to the dimension of the unipotent orbit of type $3A_1$ if $\g$
is of type $E_6$ or $E_8$ and $(3A_1)^{\prime}$ if $\g$  is of type
$E_7$. In the latter case Theorem \ref{metodo} implies that
a representative of the class of type $3A_1$ in $D_4$ is a
representative of the class of type $(3A_1)^{\prime}$. 

\medskip

\noindent
{\it Type $A_1+\tilde A_1$} ($\g$ of type $F_4$).
Let us consider the 
subgroup of $G$ of type $B_4$ generated by  $X_{\pm\alpha}$ for 
$\alpha\in\{\alpha_2+2\alpha_3+2\alpha_4,\,\alpha_1,\,\alpha_2,\,\alpha_3\}$.
By Theorem \ref{classical-unipotent}, if $\O^{\prime}$ is the
conjugacy class of 
type $A_1+\tilde A_1$ in $B_4$, then $z(\O^{\prime})$ is the longest element of the
Weyl group of $B_4$ and  coincides with the longest element of $W$.
Therefore there is a representative of the
conjugacy class of type $A_1+\tilde A_1$ in $F_4$ in $B\dot
w_0B$. We have:
$\dim\O=28=\ell(w_0)+rk(1-w_0)$.   

\medskip

\noindent
{\it Type $4A_1$} ($\g$ of type $E_7$, $E_8$). 
We observe that $\dim\O=\dim B=\ell(w_0)+rk(1-w_0)$ therefore we need to
prove that $z(\O)=w_0$.
In order to do so we shall apply Corollary \ref{reticolo}. 

Let us consider the following
subalgebras of type $D_6$ in $\g$ and the corresponding immersions of
algebraic groups: as above we realize the root systems of these subalgebras as
the sets of roots orthogonal to $\ker(1-w_i)$ where the $w_i$'s in $W$ are chosen
as follows:
\begin{itemize}
\item if $\g$ is of type $E_7$:

$w_1=s_{\beta_1}s_{\beta_2}s_{\alpha_2+\alpha_3+2\alpha_4+\alpha_5}s_{\alpha_3}s_{\alpha_2}s_{\alpha_5}=w_0s_{\alpha_7}$;

$w_2=s_{\beta_1}s_{\beta_2}s_{\alpha_2+\alpha_3+2\alpha_4+\alpha_5}s_{\alpha_3}s_{\alpha_2}s_{\alpha_7}=w_0s_{\alpha_5}$;

\item if $\g$ is of type $E_8$:

$w_1=s_{\beta_1}s_{\beta_2}s_{\beta_3}s_{\alpha_2+\alpha_3+2\alpha_4+\alpha_5}s_{\alpha_2}s_{\alpha_5}=w_0s_{\alpha_3}s_{\alpha_7}$;  

$w_2=s_{\beta_1}s_{\beta_2}s_{\beta_3}s_{\alpha_2+\alpha_3+2\alpha_4+\alpha_5}s_{\alpha_3}s_{\alpha_7}=w_0s_{\alpha_2}s_{\alpha_5}$. 
\end{itemize}

It is shown in Theorem \ref{classical-unipotent} that if $\O^{\prime}$
is the conjugacy class of
type $4A_1$ in $D_6$ then $z(\O^{\prime})$ is  the longest element of
the Weyl group 
of $D_6$ which coincides with $w_i$ in each case. The only element in
$W$ which is greater than or equal to both $w_1$ and $w_2$ is $w_0$,
hence the statement.\hfill$\Box$

\subsection{Semisimple conjugacy classes}\label{SC}
As for spherical unipotent conjugacy classes we establish a result
concerning the intersections $B^-\cap \O\cap
B\dot{w}B$, with $w\in W$, when $\O$ is a semisimple conjugacy class.

\begin{lemma}\label{manchester}
Let $t$ be a semisimple element of $G$ such that $\O_t\cap
B\dot{w}B\neq\emptyset$ for some $w\in W$. Then $B^-\cap \O_t\cap
B\dot{w}B\neq\emptyset$. 
\end{lemma}
\proof Without loss of generality we may assume that $t$ lies in
$T$. Let $g\in G$ be such that $g^{-1}tg\in B\dot{w}B$ and let
$g=u_{\sigma}\dot{\sigma}b$ be the unique decomposition of $g$ in
$U^{\sigma}\dot{\sigma}B$. 
Then 
$\dot{\sigma}^{-1}u_{\sigma}^{-1}tu_{\sigma}\dot{\sigma}$ belongs to    
$\O_t\cap B\dot{w}B\cap B^-$ since
$\dot{\sigma}^{-1}u_{\sigma}\dot{\sigma}$ lies in $U^-$. \hfill$\Box$

\subsubsection{Classical type}\label{CSS}
In this section we shall analyze the spherical semisimple conjugacy
classes of $G$ when $G$ is of classical type.
Using \cite[Remarque \S 0]{bri}
we list the spherical semisimple
conjugacy classes (up to a central element) in Table 1,
where we indicate a  
representative $g$ of each semisimple class $\O_g$, the dimension of
$\O_g$ and the structure of the Lie algebra of the  centralizer of $g$.  
By $\zeta$ we shall denote a
primitive $2n$-th root of $1$. Note that $D_1$ must be interpreted as
a $1$-dimensional torus $T_1$ wherever it occurs and that $A_0$ and
$B_0$ denote the trivial Lie algebra.

\begin{remark}\label{coniugati}{\rm It is well known that $X, Y\in Sp_{2n}$ are
conjugated in $Sp_{2n}$ if and only if they are conjugated in
$GL_{2n}$. The same holds for $X, Y$ in the orthogonal group $O_m$
(see, for example, \cite[Ex.\ 2.15 (ii)]{spst}). It follows that if $X, Y\in
SO_m$ are conjugated in $GL_m$  and $C_{O_m}(X)\not\subset SO_m$ then
$X$ and $Y$  are conjugated in $SO_m$. On the contrary, if $C_{O_m}(X)\subset SO_m$ then
the conjugacy class of $X$ in $O_m$ splits into two distinct conjugacy classes in $SO_m$ of the same dimension.}
\end{remark}

\begin{theorem}\label{semisimple}Let ${\cal O}_g$ be a
spherical semisimple conjugacy 
class of $G$.
Then $\O_g$ is \wp.
\end{theorem} 
{\center
\begin{tabular}[t]{|c|c|c|}
\hline
&$\dim\O_g$& $\mbox{Lie}(C_G(g))$\\
\hline
\hline
$A_{n-1}$&\multicolumn{2}{c|}{}\\
\hline
$\begin{array}{c}
g_k=\mbox{diag}(-I_k,\,I_{n-k})\\
k \mbox{ even and } 1\le
k\le\left[\frac{n}{2}\right]
\end{array}$
&$2k(n-k)$&${\mathbb C}+ A_{k-1}+ A_{n-k-1}$\\
\hline$\begin{array}{c}
g_{\zeta,k}=\mbox{diag}(-\zeta I_k,\,\zeta I_{n-k})
\\
k \mbox{ odd and } 1\le
k\le\left[\frac{n}{2}\right]\\
\end{array}$
&$2k(n-k)$&${\mathbb C}+ A_{k-1}+ A_{n-k-1}$\\
\hline
$B_n$&\multicolumn{2}{c|}{}\\
\hline
$\begin{array}{c}
\rho_k=\mbox{diag}(1,\, -I_k,\,I_{n-k},\,-I_k,\,I_{n-k})
\\
1\le k\le n
\end{array}$
&$2k(2n-2k+1)$&$D_k+ B_{n-k}$\\
\hline
$\begin{array}{c}
b_{\lambda}=\mbox{diag}(1,\,\lambda I_{n},\,\lambda^{-1}I_{n})
\\
\lambda\in\mathbb{C}\setminus \{0,\pm 1\}\\
\end{array}$
&$n^2+n$&${\mathbb C}+ A_{n-1}$\\
\hline
$C_n$& \multicolumn{2}{c|}{}\\
\hline
$\begin{array}{c}
\sigma_k=\mbox{diag}(-I_k,\,I_{n-k},\,-I_k,\,I_{n-k})
\\
1\le k\le\left[\frac{n}{2}\right]
\end{array}$
&$4k(n-k)$&$C_k+ C_{n-k}$\\
\hline
$\begin{array}{c}
c_{\lambda}=\mbox{diag}(\lambda, I_{n-1},\,\lambda^{-1},\,I_{n-1})
\\
\lambda\in{\mathbb C}\setminus\{0,\,\pm1\}
\end{array}$
&$4n-2$&${\mathbb C}+ C_{n-1}$\\
\hline
$c=\mbox{diag}(i\,I_{n},\,-i\,I_{n})$
&$n^2+n$&$\begin{array}{c}
{\mathbb C}+ \tilde{A}_{n-1}
 \end{array}$\\
\hline
$D_n$&\multicolumn{2}{c|}{}\\
\hline
$\begin{array}{c}
\sigma_k=\mbox{diag}(-I_k,\,I_{n-k},\,-I_k,\,I_{n-k})\\
1\le k\le\left[\frac{n}{2}\right]
\end{array}$
&$4k(n-k)$&$D_k + D_{n-k}$\\
\hline
$c=\mbox{diag}(i\,I_{n},\,-i\,I_{n})$
&$n^2-n$&${\mathbb C}+ A_{n-1}$\\
\hline
$d=\mbox{diag}(i\,I_{n-1},\,-i,\,-i\,I_{n-1},\,i)$
&$n^2-n$&${\mathbb C}+A_{n-1}$\\
\hline
\end{tabular}}

\begin{center}Table 1
\end{center}

\medskip

\Pf For each class $\O_g$ we shall exhibit an  element $w$ of the Weyl group 
such that $\dim\O_g=\ell(w)+rk(1-w)$
and a representative of $\O_g$ in $B\dot{w}B$. The proof will follow
from Lemma \ref{manchester}. 

\noindent
{\em Type $A_{n-1}$.} Let  
$t_k=\left(\begin{array}{ccc}
-\eta I_k&0&0\\
0&\eta I_{n-2k}&0\\
-2\eta Y_k&0&\eta I_k\\
\end{array}
\right)$,  where 
$\eta=\begin{cases}
\zeta&\mbox{ if } k \mbox{ is odd}\\
1&\mbox{ if } k\mbox{ is even}\\
\end{cases}$
~~~and $Y_k$ is the $k\times k$ matrix 
$\left(
\begin{picture}(30,20)
\put (0,10){0}
\put (25,10){1}
\put (0,-10){1}
\put (25,-10){0}
\put (11,1){.}
\put (15,4){.}
\put (20,7){.}
\end{picture}
\right)$.      
Then $t_k\in O_{g_k}\cap B^-$ if $k$ is even and $t_k\in\O_{g_{\zeta,k}}\cap B^-$ if
$k$ is odd. 
As in the proof of 
\cite[Theorem 3.4]{Nico2},  $t_k\in B\dot{w}_kB$
where $$w_k=(n\ n-1\ \ldots\, n-k+1\ k+1\ k+2\ \ldots n-k\ k\ \ldots 1)$$
and  $\ell(w_k)+rk(1-w_k)=2k(n-k)$.

\noindent
{\em Type $C_n$.}
Let us consider the conjugacy class $\O_{\sigma_k}$. The following
element $\dot{v}_k$  lies in  $N\cap\O_{\sigma_k}$:
$$\dot{v}_k=\left(
\begin{array}{cc|cc}
0_{2k} & 0 & S_k & 0\\
0 & I_{n-2k} & 0 & 0_{n-2k}\\
\hline
-S_k & 0 & 0_{2k} & 0\\
0 & 0_{n-2k} & 0 & I_{n-2k}
\end{array}
\right)$$
where $S_k$ is the $2k\times 2k$ matrix introduced in the proof of
Theorem \ref{classical-unipotent}. Let $v_k$ be  the image of
$\dot{v}_k$ in $W$. 
Then $$\ell(v_k)+rk(1-v_k)=4nk-4k^2=4k(n-k)=\dim{\cal
O}_{\sigma_k}.$$ 

Let us now consider the class $\O_{c_\lambda}$. Let us first assume
$n=2$. With the help of Remark \ref{bigcell}
one can check that the element
$$x=\left(\begin{array}{c|c}
A&0\\
\hline
C&D\\
\end{array}
\right)=\left(\begin{array}{cc|cc}
\lambda&0&0&0\\
0&1&0&0\\
\hline
1&1&\lambda^{-1}&0\\
\lambda&0&0&1\\
\end{array}
\right)
\in\O_{c_\lambda}\cap B^-\cap B\dot{w}_0B.$$
Let us  now suppose $n>2$. Then 
%
the element
$\left(\begin{array}{c|c}
{\begin{array}{cc}
A&\\
&I_{n-2}\end{array}} & 0_n\\
\hline
{\begin{array}{cc}
C& \\
&0_{n-2}
\end{array}}&
{\begin{array}{cc}
D & \\
&I_{n-2}
\end{array}}
\end{array}
\right)$ 
is a representative of $\O_{c_\lambda}$ lying over the following element 
$w$ of $W$:
$$w(e_i)=
\begin{cases}
-e_i&\mbox{ if }  i=1,\,2\\
e_i&\mbox{ if } i\neq1,\,2.\\
\end{cases}$$
We have: $\ell(w)+rk(1-w)=4n-2=\dim \O_{c_\lambda}$.

Let us now consider the class $\O_c$. One can check that the element $\left(\begin{array}{c|c}
0&I_n\\
\hline
-I_n&0\\
\end{array}\right)$ lies in $\O_c\cap B\dot{w}_0B$ and that $\dim\O_c=n^2+n=\ell(w_0)+rk(1-w_0)$. 

\bigskip
\noindent
{\em Type $D_n$.} Let us notice that the centralizer of $\sigma_k$ in $O_{2n}$ contains the element
$$\left(
\begin{array}{c|c}
{\begin{array}{cc}
0 & \\
 & I_{n-1} \end{array}} &{\begin{array}{cc}
1 & \\
 & 0_{n-1} \end{array}}  \\
\hline
{\begin{array}{cc}
1 & \\
 & 0_{n-1}\end{array}}&{\begin{array}{cc}
0 & \\
 & I_{n-1} \end{array}}
\end{array}
\right)$$ which does not lie in $SO_{2n}$. By Remark \ref{coniugati}, an element
 $x\in SO_{2n}$ belongs to $\O_{\sigma_k}$ if and only if it is conjugated to
$\sigma_k$ in $GL_{2n}$.

Let us consider 
the element $\tilde{w}_k$ of $W(SO_{2n})$ represented 
in $N$ by:
$$\dot{\tilde{w}}_k=\left(
\begin{array}{c|c}
{\begin{array}{cc}
0_{2k} & \\
 & I_{n-2k}
\end{array}} & {\begin{array}{cc}
I_{2k} & \\
 & 0_{n-2k}
\end{array}}\\
\hline
{\begin{array}{cc}
I_{2k} & \\
 & 0_{n-2k}
\end{array}} & {\begin{array}{cc}
0_{2k} & \\
 & I_{n-2k}
\end{array}}
\end{array}
\right).$$
Then $\dot{\tilde{w}}_k$ lies in $\O_{\sigma_k}$ and
$$\ell(\tilde{w}_k)+rk(1-\tilde{w}_k)=4k(n-k)=\dim \O_{\sigma_k}.$$
 
Let us consider the conjugacy class
$\O_c$  of $c$. In this case $C_{O_{2n}}(c)\subset SO_{2n}$ so the conjugacy class of $c$ in $O_{2n}$ splits into the conjugacy classes  $\O_c$ and  $\O_d$ 
in $SO_{2n}$. 
Let $J_k$ be the $n\times n$ matrices introduced in the proof of
Theorem \ref{classical-unipotent}.
If $n$ is even then the element $\w=\left(\begin{array}{c|c}
0&J_{n/2}\\
\hline
J_{n/2}&0\\
\end{array}\right)\in N$ is a representative of  $\O_c$. Besides,
if $w\in W$ is the image of $\w$ in $W$, then
$\ell(w)+rk(1-w)=n^2-n=\dim\O_c$. 

\noindent
If $n$ is odd then the element $\w^{\prime}=\left(\begin{array}{c|c}
{\begin{array}{cc}
0_{n-1} & \\
 & i
\end{array}}&     
J_{(n-1)/2}\\
\hline
J_{(n-1)/2} 
&
{\begin{array}{cc}
0_{n-1}& \\
 & -i
\end{array}}
\end{array}
\right)$ lies in $\O_c\cap N$ and
$\ell(w^{\prime})+rk(1-w^{\prime})=n^2-n=\dim\O_c$ where $w^{\prime}$
is the image of $\w^{\prime}$ in $W$.

Let us now consider the class $\O_d$. If $n$ is odd, then $-d\in \O_c$
so $z(\O_d)=w=z(\O_c)$. If $n$ is even, then $d=\hat\tau(c)$, where
$\hat\tau$ is the automorphism of $SO_{2n}$ introduced in the proof of
Theorem \ref{classical-unipotent}. Therefore $\hat\tau(\dot w)\in N$
is a representative of $\O_d$ and its projection $w^\tau$ is such that
$\ell(w^\tau)+rk(1-w^\tau)=n^2-n=\dim\O_d$.

\medskip
\noindent
{\em Type $B_n$.} Let $\rho_k$  with $k=1, \dots, n$ be the semisimple elements
of $SO_{2n+1}$
introduced in Table 1. The following cases need to be analysed separately:

\medskip

\noindent
{\sc Case I: $1\leq k\leq [\frac{n}{2}]$.} 
We already proved that,
under these hypotheses,
the conjugacy class
$\O_{\sigma_k}$ of $\sigma_k$ in $SO_{2n}$ contains the element
$$\dot{\tilde{w}}_k=\left( 
\begin{array}{c|c}
{\begin{array}{cc}
0 & \\
 & I_{n-2k}
\end{array}} & {\begin{array}{cc}
I_{2k} & \\
 & 0_{n-2k}
\end{array}}\\
\hline
{\begin{array}{cc}
I_{2k} & \\
 & 0_{n-2k}
\end{array}} & {\begin{array}{cc}
0 & \\
 & I_{n-2k}
\end{array}}
\end{array}
\right).$$ 
Then $\dot{v}_k=\left(\begin{array}{cc}
1 & \\
 & \dot{\tilde{w}}_k
\end{array}\right)$ lies in $\O_{\rho_k}\cap N$. 
Let $v_k$ be the element in $W(SO_{2n+1})$ represented by
 $\dot{v}_k$. Then
 $\ell(v_k)+rk(1-v_k)=2k(2n-2k+1)=\dim\O_{\rho_k}$. 

\bigskip

\noindent
{\sc Case II: $[\frac{n}{2}]<k\leq n$.}
Let us consider the following element of $N$:
$$\dot{Z}_{n-k}=\left(\begin{array}{cl}
-1 &\\
{\begin{array}{c}
\\
\\
\end{array}} &
{\begin{array}{c|c}
{\begin{array}{cc}
 ~~0_{2(n-k)+1} &\\
 & -I_{2k-n-1}
\end{array}} & {\begin{array}{cc}
 -I_{2(n-k)+1} & \\
 & 0_{2k-n-1}
\end{array}}\\
\hline
{\begin{array}{cc}
-I_{2(n-k)+1} & \\
& 0_{2k-n-1}
\end{array}} & {\begin{array}{ccc}
{\begin{array}{cc}
0_{2(n-k)+1} & \\
& -I_{2k-n-1}
\end{array}}
\end{array}}
\end{array}}
\end{array}
\right).$$ 
Since the element $\mbox{diag}(-1,I_{2n})$ belongs to the centralizer 
$C_{O_{2n+1}}(\rho_{k})$, it follows from Remark \ref{coniugati} that 
$\dot{Z}_{n-k}$ lies in $\O_{\rho_k}$. Besides,
$$\ell(Z_{n-k})+rk(1-
Z_{n-k})=2k(2n-2k+1)=\dim\,\O_{\rho_k}.$$

Finally, let $\O_{b_\lambda}$ be the conjugacy class of $b_\lambda$. Then
$$\dim\O_{b_\lambda}=n^2+n=\ell(w_0)+rk(1-w_0).$$ Let
 $$v=\left(\begin{array}{c|c|c}
1&^t\psi&0\\
\hline
0&\lambda^{-1}I_n&0\\
\hline
-\lambda\psi&\lambda\Sigma&\lambda I_n\\
\end{array}
\right)$$
where $\psi=\,^t\!(1~0~\ldots~0)$ and $\Sigma$ is the $n\times n$ matrix with diagonal 
$(-1/2,\,0,\,\ldots,\,0)$, first upper off-diagonal   
$(1,\,1,\,\ldots,\,1)$, first lower off-diagonal  
$(-1,\,-1,\,\ldots,\,-1)$ and zero elsewhere. Since the element 
$\mbox{diag}(-1,$ $I_{2n})$ belongs to the centralizer in
$O_{2n+1}$ of $b_\lambda$,  $v$ lies in $\O_{b_\lambda}$ and, by 
Remark \ref{bigcell}, lies over $w_0$.
\hfill$\Box$

\subsubsection{Exceptional Type}

In this section we shall analyze the spherical semisimple conjugacy
classes of $G$ when $G$ is of exceptional type.  
Using \cite[Remarque \S 0]{bri} we are able to list the spherical
semisimple conjugacy classes up to a central element. The results are collected in
Table 2, where we indicate a representative $g$ of each semisimple class $\O_g$, the
dimension of $\O_g$ and  the structure of the Lie algebra of the
centralizer of $g$. If $\g$ has rank $n$ we shall denote by
$\check{\omega}_i$, for $i=1,\dots, n$, the elements in ${\mathfrak h}$
defined by 
$$\langle\alpha_j, \check{\omega}_i\rangle=\delta_{ji} ~~j=1,\dots, n.$$

\begin{theorem}\label{e-semisimple}
Let ${\cal O}_g$ be a spherical semisimple  conjugacy
class. Then $\O_g$ is \wp.  
\end{theorem}
\medskip

\begin{tabular}[t]{|c|c|c|c|}
\hline
&$\dim\O_g$& $\mbox{Lie}(C_G(g))$& $w$\\
\hline
\hline
$E_6$&\multicolumn{3}{c|}{}\\
\hline
$p_1=\exp({\pi i \check{\omega}_2})$&40&$A_1+ A_5$&$w_0=s_{\beta_1}s_{\beta_2}s_{\beta_3}s_{\alpha_4}$\\
\hline
$p_2=\exp({\pi i \check{\omega}_1})$&32&${\mathbb C}+ D_5$&$s_{\beta_1}s_{\beta_2}$\\
\hline
\hline
$E_7$&\multicolumn{3}{c|}{}           \\
\hline
$q_1=\exp({\pi i\check{\omega}_2})$&70&$A_7$&$w_0$\\
\hline
$q_2=\exp(\pi i\check{\omega}_1)$&64&$A_1+D_6$&$s_{\beta_1}s_{\beta_2}s_{\alpha_2+\alpha_3+2\alpha_4+\alpha_5}s_{\alpha_3}$\\
\hline
$q_3=\exp(\pi i\check{\omega}_7)$&54&${\mathbb C}+ E_6$&$s_{\beta_1}s_{\beta_2}s_{\alpha_7}$\\
\hline
\hline
$E_8$&\multicolumn{3}{c|}{}\\
\hline
$r_1=\exp({\pi i\check{\omega}_1})$&128&$D_8$&$w_0$\\
\hline
$r_2=\exp({\pi i\check{\omega}_8})$&112&$A_1+ E_7$&$s_{\beta_1}s_{\beta_2}s_{\beta_3}s_{\alpha_7}$\\
\hline
\hline
$F_4$&\multicolumn{3}{c|}{}\\
\hline
$f_1=\exp({\pi i\check{\omega}_1})$&28&$A_1 + C_3$&$w_0$\\
\hline
$f_2=\exp({\pi i\check{\omega}_4})$&16&$B_4$&$s_{\gamma_1}$\\
\hline
\hline
$G_2$&\multicolumn{3}{c|}{}\\
\hline
$e_1=\exp({\pi i\check{\omega}_2})$&8&
$A_1+\tilde{A}_1$ & $w_0$\\
\hline
$e_2=\exp((2\pi i\check{\omega}_1)/3)$&6&$A_2$&$s_{\gamma_1}$\\
\hline
\end{tabular}
\begin{center}Table 2
\end{center}

\medskip

\noindent
\proof Let us consider the conjugacy class $\O_{p_2}$ in $E_6$ and 
 the element
$$z=\s_{\beta_1}\s_{\beta_2}x_{\beta_1}(1)x_{\beta_2}(1)
\exp(\pi i\check{\omega}_1)x_{\beta_2}(-1)x_{\beta_1}(-1)
\s_{\beta_2}^{-1}\s_{\beta_1}^{-1}.$$  
Then
$z=x_{-\beta_1}(t_1)
x_{-\beta_2}(t_2)\s_{\beta_1}\s_{\beta_2}x_{\beta_2}(1)x_{\beta_1}(1)\exp(\pi
i\check{\omega}_1)\s_{\beta_2}^{-1}\s_{\beta_1}^{-1}$ 
for some $t_1$ and $t_2$ different from zero so 
$$z=x_{-\beta_1}(2t_1)x_{-\beta_2}(2t_2)h$$
for some $h\in T$. Hence
$z$ lies in $\O_{p_2}\cap B{\s}_{\beta_1}{\s}_{\beta_2}B\cap B^-$. Besides,
$\ell(s_{\beta_1}s_{\beta_2})+rk(1-s_{\beta_1}s_{\beta_2})=
32=\dim\O_{p_2}$.

In a similar way, for the conjugacy class $\O_{q_3}$ in $E_7$, let us
consider the element 
$$y=\s_{\beta_1}\s_{\beta_2}\s_{\alpha_7}x_{\beta_1}(1)x_{\beta_2}(1)x_{\alpha_7}(1)\exp(\pi i\check{\omega}_7)x_{\alpha_7}(-1)
x_{\beta_2}(-1)x_{\beta_1}(-1)\s_{\alpha_7}^{-1} \s_{\beta_2}^{-1}\s_{\beta_1}^{-1}$$
$$=x_{-\beta_1}(t_1)x_{-\beta_2}(t_2)x_{-\alpha_7}(t_3) 
\s_{\beta_1}\s_{\beta_2}\s_{\alpha_7}x_{\alpha_7}(1) x_{\beta_2}(1)x_{\beta_1}(1)\exp(\pi i\check{\omega}_7)
\s_{\alpha_7}^{-1}
\s_{\beta_2}^{-1}\s_{\beta_1}^{-1}$$
for some $t_1$, $t_2$ and $t_3$ different from zero. Then
$$y=x_{-\beta_1}(2t_1)x_{-\beta_2}(2t_2)x_{-\alpha_7}(2t_3) h$$ for some
$h\in T$. Hence, $y$ lies in $\O_{q_3}\cap
B\s_{\beta_1}\s_{\beta_2}\s_{\alpha_7}B\cap B^-$. Besides,
$\ell(s_{\beta_1}s_{\beta_2}s_{\alpha_7})+rk(1-s_{\beta_1}s_{\beta_2}s_{\alpha_7})=54=\dim\O_{q_3}$.  

\medskip

Let now $\g$ be of type $F_4$ and let us consider the conjugacy class
$\O_{f_2}$.
 Let
us fix a short root $\gamma$ which does not belong to the root system
of $C_G(f_2)$ (which is of type $B_4$), and let $w\in W$ be such that
$w(\gamma)=-\gamma_1$. Then we have:
$$x:=\dot{w}x_{\gamma}(-1)f_2x_{\gamma}(1)\dot{w}^{-1}=x_{-\gamma_1}(t)h$$
for some $t\neq 0$ and some $h\in T$. Therefore $x$ lies in
$\O_{f_2}\cap B\dot{s}_{\gamma_1}B\cap B^-$. Since the root system of
type $F_4$ is self-dual we have:
$$\ell(s_{\gamma_1})+rk(1-s_{\gamma_1})=\ell(s_{\beta_1})+rk(1-s_{\beta_1})=16=\dim \O_{f_2}.$$

Let now $\g$ be of type $G_2$. 
Then  the element 
$$\dot{s}_{\gamma_1}^{-1}x_{\gamma_1}(-1)e_2x_{\gamma_1}(1)\dot{s}_{\gamma_1}=
\dot{s}_{\gamma_1}^{-1}x_{\gamma_1}(t)e_2\dot{s}_{\gamma_1},$$
for some $t\neq 0$, lies in $B\dot{s}_{\gamma_1}B\cap\O_{e_2}\cap B^-$.
As in type $F_4$ we have: 
$$\ell(s_{\gamma_1})+rk(1-s_{\gamma_1})=\ell(s_{\beta_1})+rk(1-s_{\beta_1})=6=\dim\O_{e_2}.$$  

For the remaining spherical semisimple conjugacy classes we shall
assume $G=G_{ad}$ and use Lemma \ref{isogeny}.  For each of these
classes $\O_g$
we shall prove the statement by exhibiting an element
$\dot{w}\in N\cap\O_g$
such that $\ell(w)+rk(1-w)=\dim\O_g$ and by using Lemma
\ref{manchester}.    
The elements $w$'s are listed in Table 2. 
Let us observe that for every element $w$ in Table 2 corresponding to
these classes
we can choose a
representative $\dot{w}\in N$ of order two in $G_{ad}$. 
For $w=w_0$, when $w_0=-1$, this fact was observed in
\cite[Lemma 2]{volk}.
In general this can be seen 
using the expression of $w$ as a  product of reflections with respect
to mutually orthogonal roots as in Table 2, and \cite[Lemma
7.2.1]{Carter1}.    
From the analysis of the
conjugacy classes of the involutions of $G_{ad}$ in \cite{Kac} (see
also \cite[\S X.5]{Helga}) we deduce
$\dim\O_{\dot{w}}\leq \ell({w}_0)+rk(1-{w}_0)$. If $w=w_0$, by
Theorem \ref{metodo},
$\dim\O_{\dot{w}_0}=\ell({w}_0)+rk(1-{w}_0)$. By \cite[\S X.5, Tables
II, III]{Helga} 
there is only one conjugacy class of involutions in $G_{ad}$ whose
dimension is equal to $\ell({w}_0)+rk(1-{w}_0)$. Therefore $\dot{w}_0$
lies in the spherical semisimple conjugacy class of maximal dimension.

Finally, we are left with the conjugacy classes $\O_{q_2}$ and
$\O_{r_2}$.  In order 
to prove that the element
$\dot{w}$ lies in the corresponding orbit $\O_g$ when $g$ is either
${q_2}$ or ${r_2}$, it is sufficient to use \cite[\S X.5, Tables
II, III]{Helga} and  estimate the dimension of the centralizer of $\w$. One can
perform this computation in the Lie algebra of $G$, namely, calculating
the dimension of $\mbox{Lie}(C_G(\dot{w}))=\{x\in\g ~|~
Ad(\dot{w})(x)=x\}$. This can be done analyzing the eigenspaces of
$Ad(\dot{w})$ in the stable subspaces of the form $\g_{\alpha}+
\g_{w(\alpha)}$, with the use of
\cite[Lemma 7.2.1]{Carter1}. \hfill$\Box$  

\subsection{The remaining conjugacy classes}
In this section we shall investigate the spherical conjugacy classes
${\cal O}_g$ of elements
$g\in G$ which are neither semisimple nor unipotent.
If the conjugacy class $\O_g$ of an element $g$ with
Jordan decomposition $su$ is spherical, then both $\O_{s}$ and
$\O_{u}$ are spherical. Indeed, if
$BC_G(g)$ is dense in $G$ then also $BC_G(s)\supset BC_G(g)$ and
$BC_G(u)\supset BC_G(g)$ are dense in $G$.  Therefore the semisimple
parts of the elements we shall consider in this section are those
occurring in \S \ref{SC}. 
Let us notice that when  
the identity component 
of the centralizer of such a semisimple element 
is not simple it is isomorphic either to an almost direct product $G_1G_2$ or to an
almost direct product
 $G_1G_2T_1$ where $T_1$ is a one-dimensional torus. 
When this is the case we will identify $G_j$ with a subgroup of
$G$ and  write a unipotent element commuting with $s$ 
as a  pair  $(u_1, u_2)$ or, equivalently, as a product $u_1u_2$ with $u_j\in G_j$ unipotent.
If the conjugacy class of
$su=su_1u_2$ is spherical, 
then the conjugacy class 
of $u_j\in G_j$  is necessarily spherical.

In the sequel we will need the following definition and results:
\begin{definition}Let $\bar G$ be a reductive connected algebraic group.
Let $H$ be a closed subgroup of $\bar G$. We say that $H=H^u K$ is a Levi decomposition
 of $H$ if $H^u$ is the unipotent radical of 
$H$ and $K$ is a maximal reductive subgroup of $H$. 
\end{definition}
In characteristic zero such a decomposition always exists.
\begin{proposition}\label{condition}\cite[Proposition I.1]{bri} Let $\bar G$ be a
 reductive connected algebraic group over an algebraically
 closed field of characteristic zero. Let $H$ be a closed 
 subgroup of $\G$ with Levi decomposition $H=H^u K$. Let
 $P$ be a parabolic subgroup of $G$ with a Levi
 decomposition $P=P^uL$ such that 
$H^u\subset P^u$ and $K\subset L$. Then 
the following conditions are equivalent:
\begin{enumerate}
\item ${\G}/H$ is spherical;
\item $K$ has an open orbit in $P^u/H^u$ and the generic $K$-stabilizer of $P^u/H^u$ is spherical in $L$.
\end{enumerate}
\end{proposition}
%
%
When $H$ is the centralizer $C_\G(u)$ of a
unipotent  element in  a semisimple algebraic group $\G$,
a construction of the subgroups $P$,  $K$ and  $L$ as in Proposition
\ref{condition} is given in 
\cite[Lemma 5.3]{elki}, using key results of \cite{spst}.
Let us recall this construction.
Let  $e$ be  the nilpotent  element of $\bar{\g}=\mbox{Lie}(\bar{G})$
corresponding to $u$ and let
$(e,h,f)$ be an $sl_2$-triple in $\gg$. 
The semisimple element $h$ determines a natural $\mathbb Z$-grading on
$\gg$ by $\gg_j:=\{z\in\gg~|~[h,z]=jz\}$. 
The  subalgebra 
$\p:=\bigoplus_{j\ge0}\gg_j$ is parabolic and 
$\p^u:=\bigoplus_{j>0}\gg_j$ is its
nilpotent radical. The subalgebra $\p$ is called the canonical
parabolic subalgebra associated to $e$
and it is independent of the choice of the $sl_2$-triple.
Let $P$ be the parabolic subgroup of $\G$ whose Lie algebra is $\p$
and let $L$ be 
the connected, reductive subgroup of $\G$  whose Lie algebra is
$\gg_0$, i.e., $L=\{g\in \G| Ad (g)h=h\}^\circ$. 
The group $P$ is called the {\em canonical parabolic} associated to
$u$ and
$P=P^uL$ is a Levi decomposition of $P$. It turns out that
$C_\G(u)\subset P$, $C_\G(u)^u\subset P^u$ 
and that $C_\G(u)=(P\cap C_\G(u))(C_G(u)\cap L)$ is a
Levi decomposition of $C_\G(u)$.
%
%

\medskip

A similar construction works in the case of non-semisimple elements:

\begin{lemma}\label{levi}Let $\G$ be a connected reductive
algebraic group with Lie algebra $\gg$, 
let $g\in \G$ be an element with Jordan decomposition $g=su$, $u\not=1$, and let 
$H=C_\G(g)$.
Then the Levi decomposition $P=P^uL$ of the  canonical parabolic $P$ 
associated to $u$ induces a Levi decomposition $H=H^uK$ of $H$ with
$K=L\cap H$. 
\end{lemma}
\proof The semisimple element  $s$ lies in $C_\G(u)$ and  
$u$ lies in $C_\G(s)^\circ$
 which is a reductive subgroup. Hence,
there exists an $sl_2$-triple 
$(e,\,h,\,f)$  of elements of 
${\mbox{Lie}}(C_\G(s)^\circ)=\{x\in\gg\ |\ Ad(s)x=x\}$ where $e$ is
the nilpotent element associated
 to $u$. It follows that $s\in C_\G(h)=\{y\in \G\ |\ Ad(y)h=h\}=\tilde L$ 
where $\tilde L^\circ=L$.
The canonical parabolic $P$ associated to
$u$  contains $H=C_\G(u)\cap C_\G(s)$.
The subgroup 
$K=L\cap H$ is reductive because it  is the centralizer
of a semisimple element  $s\in\tilde L\cap C_G(u)$  
(see \cite[Corollary 9.4]{st}). 
The subgroup $$V=P^u\cap H=C_\G(u)^u\cap C_{\G}(s)$$
is a unipotent normal subgroup of $H$. 
In order to prove that $H=KV$ is a Levi decomposition of $H$ and,
 in particular, that $H^u=V$, it is enough to show that
$H\subseteq KV$ because $K\cap V=1$ follows from the Levi decomposition of 
$C_\G(u)$.\\
Let $z\in H$. As $H\subset C_\G(u)$,
there exist unique $v\in C_\G(u)^u$ and $t\in C_\G(u)\cap L$
such that $z=vt$. Then $svs^{-1}\in V$
because $V$ is normal in $H$ and $sts^{-1}\in
C_\G(u)\cap L$  because both $t,\,s\in C_\G(u)\cap \tilde L$ and $L$
is normal in $\tilde L$. Besides, $z=szs^{-1}$. By
the uniqueness of the decomposition in $C_\G(u)$ we get necessarily
$sts^{-1}=t$ and
 $svs^{-1}=v$, i.e., 
$t\in K$ and $v\in V$.
\hfill$\Box$

\begin{corollary}\label{levic}Let $\G$ be a connected reductive
algebraic group with Lie\
 algebra $\gg$, 
let $g\in \G$ be an element with Jordan decomposition $g=su$,
 $u\neq1$,
 and let $H=C_\G(g)$.
 Then the Levi decomposition $P=P^uL$ of the  canonical parabolic $P$
associated to $u$ induces a Levi decomposition $H^\circ=H^uK^\circ$
of $H^\circ$ with $K=L\cap H$. 
\end{corollary}
\proof The corollary follows from $H^u\subseteq H^\circ$.\hfill$\Box$

\medskip

As we already observed the sphericity of $\G/H$ depends only on the
Lie algebras of $\G$ and $H$.
In particular for the analysis of the conjugacy class of an element
$g\in G$
it does not matter whether we consider $C_G(g)$ or its identity component.

\begin{remark}\label{height}{\rm Let $G_1\subset G_2$ be reductive
algebraic groups and let $u$ be a unipotent element in $G_1$. Suppose
that  the conjugacy class of  $u$ in $G_2$ is spherical. Then the
conjugacy class of $u$ in $G_1$ is spherical by \cite[Corollary 2.3,
Theorem 3.1]{pany}}. \end{remark}

Again we shall handle
the classical and the exceptional cases separately.
\subsubsection{Classical type}
In this section we shall assume that $G$ is of classical type.
\begin{proposition}\label{clsf}
Let  $g=su\in G$ with $s\neq 1$ and $u\neq1$. 
If the conjugacy class of $g$ is spherical then only the following
possibilities may occur:
\begin{itemize} 
\item $G$ is of type $C_n$ and, up to a central element,  $g=\sigma_k 
u$ with $u=X_{1,2n}$;
\item $G$ is of type $B_n$ and, up to a central
element, $g=\rho_n u$ where $u=X_{2t, 2n+1}$ with $t=1,\dots ,
\left[\frac{n}{2}\right]$. 
\end{itemize} 
\end{proposition}
\proof We shall use Proposition \ref{condition} in order to show that
if $g$ is not as in the statement, then $\O_g$ cannot be spherical.
With notation as in Lemma \ref{levi} we shall describe $K^{\circ}$ and
its action on $P^u/H^u\cong \p^u/\h^u$.    

\noindent
{\em Type $A_{n-1}$.} Since $\O_g$ is spherical
 $s$ is conjugated, up to a central element, to one  
 of the $g_k$'s or of the $g_{\zeta, k}$'s (see Table 1). We
shall show that necessarily $u=1$, leading to a contradiction. As $u$
 is a unipotent element of the centralizer of $g_k$ (resp.\  $g_{\zeta, k}$), 
it can be identified with a pair $(u_1,\,u_2)\in SL_k\times SL_{n-k}$.
It is enough to prove that
if one of the $u_j\not=1$ and the other is equal to $1$ then $\O_g$ is
 not spherical. Suppose that $u=(u_1,\,1)$ with $u_1$ spherical with
 Young diagram of shape $X_{t,\,k}$ for $1\le
 t\le\left[\frac{k}{2}\right]$.  
We have: 
\begin{itemize}
\item $\p^u/\h^u\cong
 Mat_{t,n-k}\times Mat_{n-k,t}$; 
\item $K\cong\{
 (A,\,B,\,C)\in GL_t\times GL_{k-2t}\times GL_{n-k} ~|~ \det
 A^2\det B\det C=1\}$;
\item action of $K$ on $\p^u/\h^u$: 
$$(A,\,B,\,C).(P,\,Q)=(Y_tAY_tPC^{-1},\,CQA^{-1})$$
where $Y_t$ is a symmetric $t\times t$ matrix such that $Y_t^2=1$,
depending on the choice of $u_1$. 
\end{itemize}
Since $Tr(QY_tP)$ is a non-trivial
polynomial invariant of the action of $K$ on $\p^u/\h^u$,  
$\O_g$ is not spherical. The case
$u_2\not=1$ is similar and left to the reader.

\medskip

Let now  $G$ be  orthogonal or symplectic. Then if the conjugacy class
 of $g=su_1u_2$ is spherical, then $u_1$
 and $u_2$ are either of shape
 $X_{t,m}$
 or of shape $Z_{2t,m}$, with $u_j$ of shape $Z_{2t,m}$ only
if $G=SO_m$.

\noindent
{\em Type $C_n$.} Let us distinguish the following
possibilities for $s$:

\noindent
i) $s=\sigma_k$.
If $u_1=X_{t,2k}$ with $t\geq 1$ and $u_2=1$ we have:
\begin{itemize}
\item $\p^u/\h^u\cong Mat_{2n-2k,t}$;
\item $K^\circ\cong Sp_{2n-2k}\times SO_t\times Sp_{2k-2t}$;
\item action of $K^\circ$: orthosymplectic of $Sp_{2n-2k}\times SO_t$.
\end{itemize}
If $t\ge2$ the orthosymplectic action of $Sp_{2n-2k}\times SO_t$
 on $Mat_{2n-2k,t}$ cannot have a dense orbit because it has a non-trivial
 invariant. Indeed, if $X\in Mat_{2n-2k,t}$, $E$ is the matrix of the form with respect to which  $SO_t$ is orthogonal, and if $J$ is the 
 matrix of the form with respect to which  $Sp_{2n-2k}$ is symplectic, then 
$Tr((E\,^t\!XJX)^2)$ is a non-trivial invariant for the $Sp_{2n-2k}\times SO_t$-action.
Then, if $u_1$ is of shape
$X_{t,2k}$ and $t\ge2$,  
$\O_{\sigma_ku}$ is not spherical.   By the symmetry in the roles of
$u_1$ and $u_2$ the same holds if  $u_2$ is of shape
$X_{t,2n-2k}$ with $t\geq 2$. 

If
$u_1=X_{1,2k}$ and $u_2=X_{1,2n-2k}$ we have:
\begin{itemize}
\item $\ph\cong{\mathbb C}^{2n-2k}\oplus{\mathbb C}^{2k-2}
\oplus {\mathbb C}$
\item $\K\cong Sp_{2n-2k}\times Sp_{2k-2}$
\item action of $\K$: standard of $Sp_{2n-2k}$ $\oplus$ standard of
$Sp_{2k-2}$ $\oplus$ trivial.
\end{itemize}
It is clear  that the action of 
$K^\circ$ on $\ph$ cannot have an open orbit.

\medskip

\noindent
ii) $s=c_\lambda$. Since $u\in C_G(c_{\lambda})$, 
$$u=\left(\begin{array}{cc|cc} 1&&&\\
&U_1&&U_2\\
\hline
&&1&\\
&U_3&&U_4\end{array}\right)$$ where
$\left(\begin{array}{c|c}U_1&U_2\\ \hline 
U_3&U_4\end{array}
\right)$ is a spherical unipotent element of $Sp_{2n-2}$. In
particular the Young diagram of $u$ has shape $X_{k, 2n}$ with $1\leq
k\leq n-1$. 
We have: 
\begin{itemize}
\item $\ph\cong \mathbb{C}^k\oplus
\mathbb{C}^k$;
\item $K^{\circ}\cong\mathbb{C}^*\times SO_k\times
Sp_{2n-2k-2}$;
\item  $K^{\circ}$ acts as follows:
$(a, A, B).(v,w)=(aAv, a^{-1}Aw).$
\end{itemize}
This action has never a dense orbit since the product
$(^tvEv)(^twEw)$ is invariant. 

\medskip 

\noindent
iii) $s=c$.
Then necessarily
$$u=\left(\begin{array}{c|c}
A&\\
\hline
&^tA^{-1}\end{array}\right)$$ where
$A$  is a spherical unipotent element of $SL_n$. In particular the
Young diagram of $u$ has shape $X_{2k, 2n}$ with $1\leq k\leq
[\frac{n}{2}]$.  
We have:
\begin{itemize}
\item $\ph\cong Sym_k\oplus Sym_k\oplus Mat_{k,
n-2k}\oplus Mat_{n-2k, k}$ where $Sym_k$ is the space of $k\times k$ symmetric matrices; 
\item $K^{\circ}\cong GL_k\times GL_{n-2k}$; 
\item $K^{\circ}$ acts 
as follows:
$$(A,B).(Z,L,M,N)=$$
$$=(Y_kAY_kZY_k~^t\!AY_k,~
^t\!A^{-1}LA^{-1}, ~^t\!A^{-1}MB^{-1}, BNY_k~^t\!AY_k).$$
\end{itemize}
This action has never a dense orbit
since $Tr(Y_kZY_kL)$ is a nonzero polynomial invariant.

\noindent
{\em Type $D_n$.} Let us distinguish the following
possibilities for $s$:

\noindent
i) $s=c.$
This case can be treated as for  $G$ of type $C_n$. In the
computations $Sym_k$ is replaced by $Ant_k$, the space of
skew-symmetric $k\times k$ matrices. When $k=1$ the product $MN$ is a
non-trivial invariant.

\medskip
\noindent
ii) $s=d$. If $n$ is odd the proof follows by noticing that
$\O_c=\O_{-d}$.
If $n$ is even the conjugacy class of $cu$ is spherical if and only if
the conjugacy class of $\hat\tau(cu)=d\hat\tau(u)$ is spherical. Then
the proof follows from i).

\medskip
\noindent 
iii) $s=\sigma_k$. If
$u_1=X_{2t,2k}$ and $u_2=1$ we have:
\begin{itemize}
\item $\ph\cong Mat_{2n-2k,2t}$;
\item $\K\cong SO_{2n-2k}\times Sp_{2t}\times SO_{2k-4t}$;
\item the action of $\K$ is the
 orthosymplectic of $SO_{2n-2k}\times Sp_{2t}$.
\end{itemize}
If
$u_1=Z_{2t,2k}$ and $u_2=1$ we have:
\begin{itemize}
\item $\ph\cong Mat_{2n-2k,2t}
\oplus{\mathbb C}^{2n-2k}\oplus{\mathbb C}^{2t}$;
\item $\K\cong SO_{2n-2k}\times Sp_{2t}\times SO_{2k-4t-3}$;
\item action of $\K$: orthosymplectic
of $SO_{2n-2k}\times Sp_{2t}$  $\oplus$ standard of $SO_{2n-2k}$
$\oplus$ standard of $Sp_{2t}$. 
\end{itemize}
The orthosymplectic action of $SO_{2n-2k}\times Sp_{2t}$ on $Mat_{2n-2k,2t}$
has  a non-trivial invariant, namely $Tr((J\,^t\!XEX)^2)$, unless
$t=0$ which occurs only if $u_1=Z_{2t,2k}$.

If 
$u_1=Z_{0,2k}$ the standard action of $SO_{2n-2k}$
 has no dense orbit because
if $v\in {\mathbb C}^{2n-2k}$ then 
$^tvEv$ is a non-trivial invariant. Therefore when $u_2=1$ the conjugacy
 class of $g$ is not
 spherical unless
$g=u_1$ or $g=\sigma_k$, leading to a contradiction.
 By the symmetry in the roles of $u_1$ and $u_2$ the result follows
for $G$ of type $D_n$.

\medskip
\noindent
{\em Type $B_n$.} Let us distinguish the following possibilities for $s$:

\noindent
i) $s=\rho_k$. Let $u=X_{2t,r}$, with $r=2k$ if $u_2=1$ and $r=2n-2k+1$
if $u_1=1$. We have:
\begin{itemize}
\item 
$\ph\cong Mat_{2n+1-r,2t}$;
\item $\K\cong SO_{2n+1-r}\times Sp_{2t}\times SO_{r-4t}$;
\item the action of $\K$ is the
 orthosymplectic of $SO_{2n+1-r}\times Sp_{2t}$.
\end{itemize}
Let $u=Z_{2t,r}$, with $r=2k$ if $u_2=1$ and $r=2n-2k+1$
if $u_1=1$. We have:
\begin{itemize}
\item $\ph\cong Mat_{2n+1-r,2t}\oplus{\mathbb
C}^{2n+1-r}\oplus{\mathbb C}^{2t}$;
\item $\K\cong SO_{2n+1-r}\times Sp_{2t}\times SO_{r-4t-3}$;
\item the action of $\K$ is orthosymplectic of
$ SO_{2n+1-r}\times Sp_{2t}$
 $\oplus$ standard of $ SO_{2n+1-r}$
$\oplus$ standard of $Sp_{2t}$.
\end{itemize}
In both cases, by arguments similar to the previous ones, the action of $K^\circ$ on 
$\ph$ can never have a dense orbit unless 
$g=\rho_n u$ with 
$u=X_{2t,2n+1}$. 

\medskip
\noindent
ii) $s=b_{\lambda}$.
Then, necessarily,
$$u=\left(\begin{array}{cc|c} 1&&\\
&A&\\
\hline
&&^tA^{-1}\end{array}\right)$$ where
$A$  is a spherical unipotent element of $SL_n$. In particular the
Young diagram of $u$ has shape $X_{2k, 2n+1}$ with $1\leq k\leq
[\frac{n}{2}]$. 
We have: 
\begin{itemize}
\item $\ph\cong \mathbb{C}^k\oplus
\mathbb{C}^k\oplus Ant_k\oplus Ant_k\oplus Mat_{k,
n-2k}\oplus Mat_{n-2k, k}$;
\item $K^{\circ}\cong GL_k\times GL_{n-2k}$;
\item  $K^{\circ}$ acts on $\ph$
as follows:
$$(A,B).(v,w,Z,L,M,N)=$$
$$=(^t\!A^{-1}v, Y_kAY_kw, Y_kAY_kZY_k~^t\!AY_k,~
^t\!A^{-1}LA^{-1}, ~^t\!A^{-1}MB^{-1},\, BNY_k~^t\!AY_k)$$
where $Y_k$ is as above.
\end{itemize}
This action has never a dense orbit
since $^t\!w Y_k v$ is a nonzero polynomial invariant. 
The statement now follows.\hfill$\Box$ 

\medskip

Let us now analyze the remaining possibilities.

\begin{theorem}\label{C-miste}
Let $g=su$ be an element of $G$ such that:
\begin{itemize}
\item either $G$ is of type $C_n$, $s=\sigma_k$ and $u=X_{1,2n}$;
\item or $G$ is of type $B_n$, $s=\rho_n$ and $u$ is a
spherical unipotent element associated to a Young diagram with two
columns.
\end{itemize} Then $\O_g$ is spherical and \wp.
\end{theorem}
\Pf We shall show that $\O_g$ is  \wp\ and hence spherical by  
exhibiting an element $x\in \O_g\cap B\dot{w} B\cap B^-$
for some $w$ such that $\ell(w)+rk(1-w)=\dim \O_g$.

\noindent
{\em Type $C_n$.} Let
$u=(u_1, u_2)\in C_{Sp_{2n}}(\sigma_k)\cong Sp_{2k}\times
Sp_{2n-2k}$, where $1\leq k\leq [\frac{n}{2}]$, and let us distinguish
the following cases: 

1. $u_1=1$, $u_2=X_{1,2(n-k)}$. In this case $\dim {\cal
O}_g=(4k+2)(n-k)$. 

(i) Let us assume
$k=[\frac{n}{2}]$. Then
$\dim\O_g=n^2+n=\ell(w_0)+rk(1-w_0)$.    
Let us choose the following
element $M\in B^-$:
$$M=\left(\begin{array}{c|c}
S
& ~~0~~\\
\hline
{\begin{array}{cccccc}
1 & 1 & 0 & &  \\ 
-1 & 0 & 1 & &  \\
 0 & -1 & 0 & 1  \\
 & &   -1 & \ddots &\ddots \\
 & &  & \ddots  & \ddots 
\end{array}}
 & ~~S~~
\end{array}
\right)$$
where $S=\mbox{diag}(1, -1, 1, -1, 1,  \dots)$ is a $n\times n$ matrix.
Then one can verify that $M$ lies
over $w_0$ and that $M\in\O_g$. 

\medskip
(ii)
Now let us suppose $k<[\frac{n}{2}]$. Notice that in this case $n-k\geq
k+2$. Let $i_{2k+1}$ be the following embedding of $Sp_{4k+2}$ into $Sp_{2n}$:
$$\left(\begin{array}{c|c}
A_{2k+1} & B_{2k+1}\\
\hline
C_{2k+1} & D_{2k+1}
\end{array}
\right)\stackrel{i_{2k+1}}\longmapsto \left(\begin{array}{c|c}
{\begin{array}{cc}
A_{2k+1} & \\
 & I_{n-2k-1}
\end{array}} & {\begin{array}{cc}
B_{2k+1} & \\
 & 0_{n-2k-1}
\end{array}}\\
\hline 
{\begin{array}{cc}
C_{2k+1} & \\
 & 0_{n-2k-1}
\end{array}} & {\begin{array}{cc}
D_{2k+1} & \\
 & I_{n-2k-1}
\end{array}}
\end{array}
\right).$$ 
Case (i) shows that if $G=Sp_{4k+2}$, $g^{\prime}=\sigma_ku$ where $u=(u_1,
u_2)$, $u_1=1$, $u_2=X_{1, 2(k+1)}$ then $\O_{g^{\prime}}$ contains a matrix $M\in
B^-(Sp_{4k+2})$ lying over $w_0$. In particular this implies that
$i_{2k+1}(M)$ lies in $B^-\cap B\dot{w}_{2k+1}B$ where
$$\dot{w}_{2k+1}=\left(\begin{array}{cc|cc}
0_{2k+1} & 0 & -I_{2k+1} & 0\\
0 & I_{n-2k-1} & 0 & 0_{n-2k-1}\\
\hline
I_{2k+1} & 0 & 0_{2k+1} & 0\\
0 & 0_{n-2k-1} & 0 & I_{n-2k-1}
\end{array}
\right).$$
The thesis follows by noticing that $i_{2k+1}(M)$ belongs to $\O_g$ and that
$$\ell(w_{2k+1})+rk(1-w_{2k+1})=(4k+2)(n-k)=\dim {\cal O}_g.$$ 

2. $u_1=X_{1, 2k}$, $u_2=1$. In this case $\dim {\cal
O}_g=2k(2n-2k+1)$.

(i) Let us first suppose that $n$ is even and let  $k=\frac{n}{2}$ so that $\dim \O_g=n^2+n$. 
  Let us
 choose the following
element $\bar{M}\in B^-\cap \O_g$:
$$\bar{M}=\left(\begin{array}{c|c}
\bar{S} & ~~0~~\\
\hline
{\begin{array}{ccccccc}
1 & 1 & & & & &  \\
-1 & 0 & 1 & & &\mbox{\Large{0}} &  \\
  & -1 & 0 & 1 & & &  \\
 & &  \ddots & \ddots & \ddots & & \\
 & &   &\ddots  & \ddots & 1 &\\
 &\mbox{\Large{0}} & & & -1 & 0 & -1\\
 & & & & & 1 & 0 
\end{array}}
 & ~~\bar{S}~~
\end{array}
\right)$$
where $\bar{S}=\mbox{diag}(-1, 1, -1, 1,  \dots)$ is an $n\times n$ matrix.  Since 
$\dim\O_g=\ell(w_0)+rk(1-w_0)$, it is enough to show that $\bar{M}$ lies
over $w_0$ and this follows, using Remark \ref{bigcell}, from a straightforward
calculation.

(ii) Now let us suppose $k<\frac{n}{2}$. Case (i) shows that if $G=Sp_{4k}$, $g^{\prime}=\sigma_ku$
where $u=(u_1, u_2)$, $u_1=X_{1, 2k}$, $u_2=1$, then $\O_{g^{\prime}}$ contains a matrix 
$\bar{M}\in B^-(Sp_{4k})$ lying over $w_0$. Using the embedding $i_{2k}$ of
$Sp_{4k}$ into $Sp_{2n}$, it is immediate to see that $i_{2k}(\bar{M})$ belongs to 
$B^-\cap B\dot{\sigma}B$ where
$$\dot{\sigma}=\left(\begin{array}{cc|cc}
0_{2k} & 0 & -I_{2k} & 0\\
0 & I_{n-2k} & 0 & 0_{n-2k}\\
\hline
I_{2k} & 0 & 0_{2k} & 0\\
0 & 0_{n-2k} & 0 & I_{n-2k}
\end{array}
\right).$$
Finally, let us notice that $i_{2k}(\bar{M})$ is conjugated to $g$ and that
if $\sigma$ is the projection of $\dot{\sigma}$ in $W$, then
$\ell(\sigma)+rk(1-\sigma)=4kn-4k^2+2k=\dim{\cal O}_g.$

\medskip
\noindent
{\em Type $B_n$.} Let $g=\rho_n u$ where $u=(u_1,1)$, $u_ 1$ is of shape
$X_{2k,2n}$ and 
$1\leq k\leq \left[n/2\right]$. If $k<\frac{n}{2}$ the class $\O_g$ is completely
determined by the diagram $X_{2k,2n}$. If $n$ is even, let $t_n\in
SO_{2n+1}$ be a representative of $s_{\alpha_n}\in 
W(SO_{2n+1})$. Conjugation by $t_n$ fixes $\rho_n$ and  induces the
automorphism $\hat{\tau}$ of $SO_{2n}$ (see the proof of
Theorem \ref{classical-unipotent}). 
Therefore, if $u_1$ and $u_1'$ are
representatives
 of the two distinct unipotent conjugacy classes of $SO_{2n}$
 associated to $X_{n,2n}$ and if $u=(u_1,1)$ and $u'=(u_1',1)$, then
$u'\in\O_u$ and $\rho_n u'\in\O_{\rho_n u}$. Therefore also for
 $k=\frac{n}{2}$ the class $\O_g$ is completely determined by
 the diagram $X_{2k,2n}$. Thus let us denote by $\O_{k}$ the conjugacy
class 
of $g=\rho_n u$ with $u_ 1$ of shape
$X_{2k,2n}$. Then $\dim\O_k=4nk-4k^2+2n-2k$.  
Let us first assume that  $k$ is maximal, i.e., $k=k_{\max}=[\frac{n}{2}]$.
Then $\dim \O_{k_{\max}}=n^2+n=\dim B(SO_{2n+1})$.\\
 Let 
$g_n=\left(\begin{array}{c|c|c}
1&\,^t\!\psi&0\\
\hline
0&-I_n&0\\
\hline
\psi&\Sigma&-I_n\\
\end{array}\right)$ where $\psi=\,^t\!(1 ~0 ~\cdots ~0)$, $\Sigma$ is the 
$n\times n$ matrix with diagonal 
$(1/2,0,\ldots,\,0)$, first upper off-diagonal $(1,\,1,\,\ldots,\,1)$, first lower off-diagonal $(-1,\,-1,\,\ldots,\,-1)$ and $0$ elsewhere. By Remark \ref{bigcell}, $g_{n}$ lies over $w_0$.\\ As $\mbox{diag}(-1,1,\,\ldots,\,1)\in C_{O_{2n+1}}(\rho_n u)$, it follows from Remark \ref{coniugati} that $g_n$
 belongs to $\O_{k_{\max}}$, so the assertion is proved for $k$ maximal.

Let us now assume that $2k< n-1 $, i.e., that there are strictly
 more than two rows
with one box in $X_{2k,2n}$.
We consider the following
embedding of $SO_{4k+3}\times SO_{2n-4k-2}$ in $SO_{2n+1}$:
$$ 
\left(\left(\begin{array}{c|c|c}
a&^t\!\alpha&^t\!\beta\\
\hline
\gamma&A&B\\
\hline
\delta&C&D\\
\end{array}\right),\,\left(\begin{array}{c|c}
A'&B'\\
\hline
C'&D'\\
\end{array}\right)\right)\mapsto
\left(\begin{array}{c|cc|cc}
a&^t\!\alpha&0&^t\!\beta&0\\
\hline
\gamma&A&0&B&0\\
0&0&A'&0&B'\\
\hline
\delta&C&0&D&0\\
0&0&C'&0&D'\\
 \end{array}\right).
$$
Let $g_{2k+1}$ be the representative of the conjugacy class
 of $\O_{k_{\max}}$ in $SO_{4k+3}$.
One can check that the embedded image of  $(g_{2k+1},-1)$
 is a representative of $\O_k$ in $B^-(SO_{2n+1})$ and that it lies over 
$\omega_k=\left(\begin{array}{c|c}
-I_{2k+1}&0\\
\hline
0&I_{n-2k-1}\\
\end{array}\right)\in W(SO_{2n+1})$.\\
 As $rk(1-\omega_k)+\ell(\omega_k)=(2k+1)^2+(2k+1)+2(n-2k-1)(2k+1)=\dim \O_k$, we have the statement for $k=1,\,\ldots,\,\left[n/2\right]$.\hfill$\Box$

\subsubsection{Exceptional type}
In this section we shall assume that $G$ is of exceptional type.
We already recalled that if the conjugacy class $\O_g$ of an element $g$ with
Jordan decomposition $su$ is spherical, then both $\O_{s}$ and
$\O_{u}$ are spherical. Besides, as $\O_g$ is spherical,
$\dim\O_g\leq\dim B$. Therefore a dimensional argument rules out all the
possibilities except the following: 

\begin{itemize}
\item $g=p_1x_{\beta_1}(1)$ if $\g$ is of type $E_6$; 
\item $g=q_2x_{\beta_1}(1)$ if $\g$ is of type $E_7$;
\item $g=r_2x_{\beta_1}(1)$ if $\g$ is of type $E_8$;
\item $g=f_2x_{\beta_1}(1)$ if $\g$ is of type $F_4$.
\end{itemize}

The following result excludes the first three cases:

\begin{proposition} If $\g$ is of type $E_6,\,E_7$ or $E_8$ any
spherical conjugacy class of $G$ is either semisimple or unipotent.
\end{proposition}
\proof By the discussion above it is enough to prove that
the class of  $sx_{\beta_1}(1)$, with $s=p_1,\,q_2,\,r_2$, is not 
spherical. Let $H$
be the centralizer of $sx_{\beta_1}(1)$ in $G$. We shall 
use the same notation as in 
Lemma \ref{levi}. Let $S$ be a stabilizer in general position for the
action of $K$ on $\l/\k$, where $\l=\mbox{Lie}(L)$ and
$\k=\mbox{Lie}(K)$.
Let $c_M(X)$ denote the complexity of the
action  of a reductive algebraic group $M$, with Borel subgroup $B_M$,
on the variety $X$, i.e., 
$c_M(X)=\min_{x\in X}{\mbox{codim}}B_M.x$. Then, by \cite[Theorem 1.2
(i)]{pany}, 
\begin{equation}c_G(G/H)=c_L(L/K)+c_S(\ph)\label{complexity}.
\end{equation}
We see that in all cases
$\l=\k\oplus \mathbb{C}h_{\beta_1}$  so that $c_L(L/K)=0$ and $S=K$. In
particular, if $\g$ is of type $E_6,\,E_7,\,E_8$  then  $K$ is of type
$A_5,\,D_6,\,E_7$, respectively. By \cite[Th\'eor\`eme 1.4]{brion2} 
(see also \cite[Theorem 1.4]{leahy}) $D_6$ and $E_7$
have no linear multiplicity free representations, hence $E_7$ and $E_8$
have no spherical exceptional conjugacy classes which are neither
semisimple nor unipotent. 

As far as $E_6$ is concerned, one can check that 
$$\p^u\simeq
\bigoplus_{\alpha>0,\atop \alpha\not\perp\beta_1}\g_{\alpha};\hskip1cm
\h^u=\g_{\beta_1},$$ therefore $\dim(\p^u/\h^u)=20$. By
\cite[Th\'eor\`eme 1.4]{brion2} there are no
multiplicity free representations  of a group of type $A_5$ on a
vector space of dimension  $20$, hence the statement.\hfill$\Box$

\begin{theorem}\label{f4} Let $\g$ be of type $F_4$ and let $\O$ be
the conjugacy class of $f_2x_{\beta_1}$. Then $\O$ is spherical and
\wp.
\end{theorem}
\Pf We have: $\dim\O=28=\ell(w_0)+rk(1-w_0)$. We shall show that $z(\O)=w_0$ which implies, by Theorem
\ref{metodo}, that $\O$ is spherical.
 
The element $f_2$ lies in $T\subset C_G(f_1)$. Besides, $C=C_G(f_1)$
 is the subgroup of $G$ of type $C_3\times A_1$ with
simple roots $\{\alpha_2,\,\alpha_3,\,\alpha_4\}$ and $\beta_1$. Since
 $(f_2)^2=1$, it follows that 
$f_2$ is of the form $(s,t)\in C_3\times A_1$  with $t$ central and $s^2=1$.
Hence,  
$f_2$ is conjugated (up to a central element) in $C$ to an element of the form
 $(\sigma_1,t)$. By  Theorem \ref{e-semisimple} 
$f_2$ is conjugated, up to a central element in $C$, by an element in the component
 of type $C_3$, to 
$\dot{s}_{\alpha_4}\dot{s}_{\alpha_2+2\alpha_3+\alpha_4}h$ for some $h\in
 T$.  Hence  
$f_2x_{-\beta_1}(1)$ is conjugated to 
$$\dot{s}_{\alpha_4}\dot{s}_{\alpha_2+2\alpha_3+\alpha_4}hx_{-\beta_1}(1)\in
B    
\dot{s}_{\alpha_4}\dot{s}_{\alpha_2+2\alpha_3+\alpha_4}\dot{s}_{\beta_1}B=B
\dot{w}_0\dot{s}_{\alpha_2}B$$ for some $h\in T$.

On the other hand, 
the involution
$\rho_4=h_{\alpha_2}(-1)h_{\alpha_2+2\alpha_3+2\alpha_4}(-1)$ (notation as in
\cite[Lemma 28]{yale}) is conjugated
to $f_2$ since  its centralizer is the subgroup of type $B_4$
with simple roots $\{\alpha_2+2\alpha_3, \alpha_1,
\alpha_2, \alpha_3+\alpha_4\}$. Therefore the element $\rho_4
x_{\beta_1}(1)\in C_G(\rho_4)$ is a representative of the class
$\O$. By Theorem \ref{C-miste}  there exists a representative of the
conjugacy class $\O_{\rho_4 x_{\beta_1}(1)}$ in $C_G(\rho_4)$ lying
over the element $w_0s_{\alpha_3+\alpha_4}$. By Corollary
\ref{reticolo} $z(\O)=w_0$. 
Let us finally show that $\O\cap B\dot{w}_0B\cap
B^-\neq\emptyset$. Let $g\in G$ be such that
$g^{-1}f_2x_{\beta_1}(1)g\in B\dot{w}_0B$ and let $g=u_\sigma\dot\sigma b$ be 
its unique decomposition in $U^\sigma\dot\sigma B$. Then
$g^{-1}f_2x_{\beta_1}(1)g$ lies in $B\dot{w}_0B$ if and only if 
$$\dot\sigma^{-1}u_{\sigma}^{-1}f_2x_{\beta_1}(1)u_\sigma\dot\sigma
=\dot\sigma^{-1}u_{\sigma}^{-1}f_2 u_\sigma\dot\sigma 
x_{\sigma^{-1}(\beta_1)}(t),$$ with $t\in {\mathbb C}^*$, lies in  $B\dot{w}_0B$.
Notice that $u_{\sigma}$ and $x_{\beta_1}(1)$ commute because
$\beta_1$ is the highest root of $\g$. The root $\sigma^{-1}(\beta_1)$ is
 negative otherwise $z(\O_{f_2})$ would be $w_0$ which is impossible
by Theorem \ref{metodo}. Then, as in Lemma \ref{manchester}, 
$\dot\sigma^{-1}u_{\sigma}^{-1}f_2 x_{\beta_1}(1) u_\sigma\dot\sigma$
lies in $B\dot{w}_0B\cap B^-$. \hfill$\Box$

\subsection{Classification and remarks}\label{conclusions}
The results of the previous sections can be summarized in the
following theorem:
\begin{theorem}\label{generale}
A conjugacy class $\O$ is spherical if and only if it is \wp.
\end{theorem} 
In fact our results lead also to the following characterization of
spherical conjugacy classes:
\begin{theorem}
Let $\O$ be a conjugacy class in   $G$, $z=z(\O)$.  Then
$\O$ is spherical
if and only if
$\dim \O = \ell(z)+rk(1-z)$.
\end{theorem}

\begin{corollary}
Let $\O$ be a spherical conjugacy class of $G$ and let $z=z(\O)$.
Let $x\in\O$ be an element such that $B.x$ is dense in
$\O$. Then $B.x=B^z.x=\O\cap B\dot{z}B$.
\end{corollary}
\proof  Theorem \ref{generale} and
Theorem \ref{metodo} show  that
if $y$ lies in $\O\cap B\dot{z}B$ then $B.y$ is dense in $\O$. It
follows that $y$ belongs to $B.x$ hence $B.x=\O\cap
B\dot{z}B$. Besides, $U^z.x=U.x$ since they are irreducible, closed
and have the same dimension. Therefore $B^z.x=TU^z.x=TU.x=B.x$.\hfill$\Box$

\bigskip

Let us introduce the map
$$\tau: ~\{\mbox{Spherical conjugacy classes of $G$}\} \longrightarrow W$$
$$\O \longmapsto z(\O)$$
and let us analyze some of its properties. A description of the image
of $\tau$ is given in Tables 2,3,4,5.
In the tables we use the notation introduced in \S \ref{E-U}.
When $G$ is of type $B$ (resp.\ $D$) the root system orthogonal to $\beta_1$
is no longer irreducible: it consists of three components of type
$A_1$ if $G$ is of type $D_4$, and of one component of type $A_1$ and
one component of type $B$ (resp.\ $D$) in the other cases. When $G$ is
of type $D_4$ we shall define $\mu_1=\alpha_1$.
When $G$ is not of type $D_4$ we shall denote by $\mu_1$ the
positive root of
the component of type $A_1$ and by  $\nu_1$ the highest root of the component of
type $B$ (resp.\ $D$). 
Inductively, for $r>1$, we shall denote by $\mu_r$, the
positive root of the component of type $A_1$ and  by $\nu_r$ the
highest root of the component of type $B$ (resp.\ $D$) of the 
root system orthogonal to $\beta_{1}, \mu_j, \nu_j$ for every
$j=1,\dots, r-1$.

In a similar way when $G$ is of type $C$ the root system orthogonal to
$\gamma_1$  consists of one component of type $A_1$ and
one component of type $C$. We shall denote by $\gamma_2^{\prime}$ the
highest short root of the component of type $C$. Inductively, for $r>1$, we shall denote by $\gamma^{\prime}_r$, the
highest short root of the component of type $C$ of the 
root system orthogonal to $\gamma_{1}, \gamma^{\prime}_j$ for every
$j=1,\dots, r-1$.

\begin{remark}\label{involution} {\rm 
We note that if $\O$ is a spherical conjugacy
class then $z(\O)$ is an involution. The reason for this
is that if $G$ is of adjoint type, then each spherical conjugacy class
$\O$ coincides with its inverse. For unipotent classes this follows
from \cite[Lemma 1.16]{Cuno}, \cite[Lemma 2.3]{Cdue}. For the
semisimple classes in 
almost all cases we are dealing with involutions in $G$. In the remaining
cases we always have $w_0=-1$, and in this case every semisimple element is
conjugate to its inverse. Finally, for the classes $\O_{su}$ where
$s\neq 1\neq u$ the
result follows from the fact that $s$ is an involution and that $u$ is
conjugate to its inverse in $C_G(s)$.}
\end{remark}

\begin{remark}
{\rm Let $\pi_1: G \longrightarrow G/U$ and $\pi_2: G/U \longrightarrow
G/B$ be
canonical projections.
Let $B$ act on $G$ by conjugation, on $G/B$ by left multiplication and on $G/U$ as follows:
$$b(gU)=bgb^{-1}U.$$
Then $\pi_1$ and $\pi_2$ are $B$-equivariant maps. In particular
$\pi_2\circ\pi_1$ maps every $B$-orbit of $G$ to a $B$-orbit of $G/B$,
i.e., a Schubert cell $C_{\sigma}=B\dot{\sigma}B/B$, for some $\sigma\in W$.

Let $\O$ be a spherical conjugacy class and let $z=z(\O)$. Let $B.x$
be the dense $B$-orbit in $\O$. Then $\dim\O=\dim
B.x=\ell(z)+rk(1-z)$. Besides $\pi_2\circ \pi_1(B.x)=C_z$ and by
\cite[Proposition 16.4]{DC-K-P3} $\dim\pi_1(\O)=\ell(z)+rk(1-z)$. It
follows that the map $\rho=\pi_1|_{\O}: \O \rightarrow G/U$
has finite fibers.
We think that the map $\rho$ could give a relation between $\O$
and the symplective leaves of $B^-$ coming from the quantization of $B^-$
(see \cite{DC-P}).}
\end{remark}

\begin{center}
\begin{tabular}[t]{|c||c|c||c|c|}
\hline
$\g$& $\O$ & $z(\O)$ & $\O$ & $z(\O)$\\
\hline\hline
$\!\!A_{n-1}\!\!$ & $X_{k,n}$ & $s_{\beta_1}\dots s_{\beta_k}$ &\multicolumn{2}{c|}{}\\
\hline
$B_n$& $X_{2k,2n+1}$ & $s_{\beta_1}s_{\nu_1}\dots s_{\nu_{k-1}}$
& $\!\!\!\!\begin{array}{c} Z_{2k,2n+1} \\  
k<\frac{(n-1)}{2}\end{array}\!\!\!\!$ & $s_{\gamma_1}\dots s_{\gamma_{2k+2}}$\\
\cline{2-5}
& $\!\!Z_{n-1, 2n+1}\!\!$ & $w_0$ &\multicolumn{2}{c|}{}\\
\hline
$C_n$& $X_{2k,2n}$ & $s_{\beta_1}\dots s_{\beta_k}$&\multicolumn{2}{c|}{}\\
\hline
$D_n$& $\!\!{\begin{array}{c}
X_{2k,2n}\\
k<\frac{n}{2}
\end{array}}\!\!$ &  $s_{\beta_1}s_{\nu_1}\dots s_{\nu_{k-1}}$ &
$Z_{2k,2n}$ &  $\!s_{\beta_1}s_{\mu_1}s_{\nu_1}s_{\mu_2}\dots s_{\nu_k}s_{\mu_{k+1}}\!$\\
\cline{2-5}
& $X_{n,2n}$ & $s_{\beta_1}s_{\nu_1}\dots s_{\nu_{n/2-2}}s_{\alpha_n}$
& $X^{\prime}_{n,2n}$ & $s_{\beta_1}s_{\nu_1}\dots s_{\nu_{n/2-2}}s_{\alpha_{n-1}}$\\
\hline
$E_6$&$A_1$ & $s_{\beta_1}$
&$2A_1$ & $s_{\beta_1}s_{\beta_2}$\\
\cline{2-5}
&$3A_1$ & $w_0$&\multicolumn{2}{c|}{}\\ 
\hline
$E_7$&$A_1$ & $s_{\beta_1}$
&$2A_1$ & $s_{\beta_1}s_{\beta_2}$\\
\cline{2-5}
&$(3A_1)^{\prime}$ & $\!s_{\beta_1}s_{\beta_2}s_{\alpha_2+\alpha_3+2\alpha_4+\alpha_5}s_{\alpha_3}\!$
&$(3A_1)^{\prime\prime}$ & $s_{\beta_1}s_{\beta_2}s_{\alpha_7}$\\
\cline{2-5}
&$4A_1$ & $w_0$&\multicolumn{2}{c|}{}\\ 
\hline
$E_8$&$A_1$ & $s_{\beta_1}$
&$2A_1$ & $s_{\beta_1}s_{\beta_2}$\\
\cline{2-5}
&$3A_1$ & $s_{\beta_1}s_{\beta_2}s_{\beta_3}s_{\alpha_7}$
&$4A_1$ & $w_0$\\ 
\hline
$F_4$&$A_1$ & $s_{\beta_1}$
&$\tilde{A}_1$ & $s_{\beta_1}s_{\beta_2}$\\
\cline{2-5}
&$A_1+\tilde{A}_1$ & $w_0$&\multicolumn{2}{c|}{}\\
\hline
$G_2$&$A_1$ & $s_{\beta_1}$
&$\tilde{A}_1$ & $w_0$\\
\hline
\end{tabular}
$$\mbox{Table 3: Unipotent
spherical conjugacy classes}$$
\end{center}

\bigskip

\begin{center}
\begin{tabular}[t]{|c||c|c|}
\hline
$\g$&$\O_{su}$& $z(\O_{su})$ \\
\hline
\hline
$B_n$& $s=\rho_n$, $u=(X_{2k,2n},1)$, $k=\left[n/2\right]$ & $w_0$\\
\cline{2-3}
&$s=\rho_n$, $u=(X_{2k,2n},1)$, $k<\left[n/2\right]$& $s_{\gamma_1}\dots
s_{\gamma_{2k+1}}$\\
\hline
$C_n$& $s=\sigma_k$, $u=(1, X_{1,2n-2k})$ & $s_{\beta_1}\dots s_{\beta_{2k+1}}$ \\
\cline{2-3}
& $s=\sigma_k$, $u=(X_{1,2k},1)$ & $s_{\beta_1}\dots s_{\beta_{2k}}$ \\
\hline
$F_4$&$f_2x_{\beta_1}(1)$&$w_0$\\
\hline
\end{tabular}
\end{center}

\begin{center}Table 4: Spherical conjugacy classes which are neither
semisimple nor unipotent
\end{center}

\medskip

\begin{center}
\begin{tabular}[t]{|c|c|c|}
\hline
$\g$&$\O$& $z(\O)$ \\
\hline
\hline
$A_{n-1}$& $\O_{g_k}$ & $s_{\beta_1}\dots s_{\beta_k}$\\
\cline{2-3}
&$\O_{g_{\zeta,k}}$ & $s_{\beta_1}\dots s_{\beta_k}$\\ 
\hline
$B_n$&$\begin{array}{c}
\O_{\rho_k}
\\
1\le k\le[n/2]
\end{array}$
& $s_{\gamma_1}\dots s_{\gamma_{2k}}$ \\
\cline{2-3}
&$\begin{array}{c}
\O_{\rho_k}
\\
\left[n/2\right] < k\le n
\end{array}$ &  $s_{\gamma_1}\dots s_{\gamma_{2(n-k)+1}}$ \\
\cline{2-3}
&$\O_{b_{\lambda}}$&$w_0$\\
\hline
$C_n$&$\O_{\sigma_k}$ & $s_{\gamma_1}s_{\gamma_2^{\prime}}\dots s_{\gamma_k^{\prime}}$\\
\cline{2-3}
&$\O_{c_{\lambda}}$ & $s_{\beta_1}s_{\beta_2}$ \\
\cline{2-3}
&$\O_c$ & $w_0$\\
\hline
$D_n$& $\O_c$ ~($n$ even) & $s_{\beta_1}s_{\nu_1}\dots
s_{\nu_{n/2-2}}s_{\alpha_n}$ \\
\cline{2-3}
 & $\O_d$ ~($n$ even) & $s_{\beta_1}s_{\nu_1}\dots
s_{\nu_{n/2-2}}s_{\alpha_{n-1}}$\\
\cline{2-3}
 & $\O_c$ ~($n$ odd) & $s_{\beta_1}s_{\nu_1}\dots s_{\nu_{(n-3)/2}}$ \\
\cline{2-3}
 & $\O_d$ ~($n$ odd) & $s_{\beta_1}s_{\nu_1}\dots s_{\nu_{(n-3)/2}}$\\
\cline{2-3}
& ${\begin{array}{c}
\O_{\sigma_k}\\
k<n/2
\end{array}}$ & $s_{\beta_1}s_{\mu_1}s_{\nu_1}s_{\mu_2}\dots
s_{\mu_{k-1}}s_{\nu_{k-1}}s_{\mu_k}$\\
\cline{2-3}
&$\O_{\sigma_{n/2}}$& $s_{\beta_1}s_{\mu_1}s_{\nu_1}s_{\mu_2}\dots
s_{\nu_{n/2-2}}s_{\mu_{n/2-1}}s_{\alpha_{n-1}}s_{\alpha_{n}}$\\
\hline
\end{tabular}
\end{center}

\begin{center}Table 5: Spherical semisimple conjugacy classes, $\g$ of
classical type
\end{center}

\medskip\medskip

\noindent
Let us recall that also $w_0$ can be decomposed as a product of
mutually orthogonal roots. 
\begin{remark}\label{5.2}{\rm We recall that for a $B$-variety $X$
 the following objects are defined:
$${\cal P}=\{f\in k(X)\setminus\{0\}\ |\  b.f=\lambda_f(b)f,\ \forall\; b\in B\}$$
where $\lambda_f\in\chi(B)$, the character group of $B$;
$$
\begin{aligned}
\psi\colon \cal P&\to \chi(B)\\
f&\mapsto \lambda_f;
\end{aligned}
$$
$$\Gamma(X)=\psi({\cal P});$$
$$r(X):=rank(\Gamma(X));$$
$$u(X)=\max_{x\in X}\dim U.x.$$

Here we note that if $\O$ is a spherical conjugacy class  then
$r(\O)=rk (1-z(\O))$ and $u(\O)=\ell(z(\O))$. Indeed, this follows from 
\cite[Corollary 1, Corollary  2(ii)]{pany3}, Theorem \ref{metodo} and 
\cite[Lemma 2.1]{knop}. }
\end{remark}

\begin{remark}{\rm Let us recall that a nilpotent orbit $\mathfrak{O}$ in $\g$ is
called a {\em model orbit} if ${\mathbb C}[{\mathfrak{O}}]$ consists exactly of
the self-dual representations of $G$ with highest weights in the root
lattice, each occurring once (see \cite[p.\ 229]{MG}). In this case the
corresponding unipotent conjugacy class  $\O$ in $G$
is spherical (\cite{brion3}, \cite{vinberg}) and, by Remark \ref{5.2}, 
$rk(1-z(\O))=rk(1-w_0)$. It follows from the proofs of
Theorems \ref{classical-unipotent}  and  \ref{exceptional-unipotent}
that $z(\O)=w_0$ 
(cf.\ \cite[Table 4]{MG2} and \cite{vogan}).}\end{remark}

\section{The proof of the DKP-conjecture}\label{conjecture}
In this section we prove the De Concini-Kac-Procesi conjecture for
representations corresponding to spherical conjugacy classes. 

Let $\ell$ be a positive odd integer greater than one.
We
will assume that $\ell$ is a good integer, i.e., that $\ell$ is
coprime with the bad primes (for the definition of the bad primes see
\cite{bourbaki}) and that $G$ is simply connected.
\subsection{Strategy of the proof} 
Let
$\varepsilon$ be a primitive $\ell$-th root of unity and
let ${\cal U}_{\varepsilon}({\mathfrak g})$ be the simply connected 
quantum group associated to ${\mathfrak g}$ as defined
in \cite{DC-K-P}, with generators $E_i, F_i, K_{\beta}$ with $\beta$
in the weight lattice $P$ and $i=1,\dots, n$. For our purposes it is convenient to introduce
the subalgebra $B_{\varepsilon}$ of  ${\cal
U}_{\varepsilon}({\mathfrak g})$ generated by $E_1, \dots, E_n$ and
$K_{\beta}$  with $\beta\in P$. The representation theory of this
algebra has been deeply investigated in \cite{DC-P}, where
$B_{\varepsilon}\cong F_{\varepsilon}[B^-]$.

The centre of ${\cal U}_{\varepsilon}({\mathfrak g})$
contains a proper, finitely generated subalgebra $Z_0$ such that  
${\cal U}_{\varepsilon}({\mathfrak g})$ is a finite $Z_0$-module (in
particular it follows that every irreducible ${\cal
U}_{\varepsilon}({\mathfrak g})$-module has finite dimension).

For any  associative algebra $A$ let us denote by $\Spec A$ the set of the equivalence classes
of the irreducible representations of $A$. It is
worth noticing that $\Spec Z_0=\{(t^{-1}u^-, tu) ~|~ u\in U, ~t\in T,
~u^-\in U^-\}$ (\cite[\S 4.4]{DC-K-P2}). In \cite{DC-K-P} the map $\pi: \Spec Z_0 \longrightarrow G$, $\pi(t^{-1}u^-, tu)=(u^-)^{-1}t^2u$,
which is an unramified covering of the big cell $\Omega=B^-B$ of $G$, is considered. Let $\varphi$ be
the map obtained by composing $\pi$ with the central character $\chi: \Spec \Ue
\longrightarrow \Spec Z_0$. It follows that for every $g\in \Omega$
one can define a certain finite-dimensional quotient ${\cal U}^g$ of
$\Ue$ such that if $g=\varphi(V)$ then $V$ is an ${\cal U}^g$-module.
 
In \cite[\S 6.1, Proposition (a)]{DC-K-P} the following crucial result is established:
\begin{equation}
\begin{array}{l}
\mbox{if}~
g, h\in\Omega ~\mbox{are conjugated in $G$ up to a central element then}\\
~{\cal U}^g ~\mbox{and} ~{\cal U}^h
~\mbox{are
isomorphic.}
\end{array}
\label{isomorfismo}
\end{equation}

In \cite[\S 6.8]{DC-K-P} the following conjecture is formulated:

\bigskip
\noindent
{\bf Conjecture}
{\it If $\sigma\in$  $\Spec{\cal U}_{\varepsilon}(\g)$ is an irreducible representation of ${\cal U}_{\varepsilon}(\g)$
on a vector space $V$ such that $\varphi (\sigma)$ belongs to a conjugacy class
${\cal O}_V$ in $G$, then $\text{dim}(V)$ is divisible by $\ell^{\frac{1}{2}\,{\text {\small dim}}\,{\cal{O}}_V}$.}

\bigskip
The De Concini-Kac-Procesi conjecture has been proved in
the following cases:
\begin{itemize}
\item[(i)] $\O$ is a regular conjugacy class (\cite[Theorem 5.1]{DC-K-P2});
\item[(ii)] $G$ is of type $A_n$ and $\ell=p$ is a prime (\cite{Nico1});
\item[(iii)] $G$ is of type $A_n$, $\ell=p^k$ and $\O$ is a
subregular unipotent conjugacy class (\cite{Nico2});
\item[(iv)] $G$ is of type $A_n$ and $\O$ is a spherical unipotent
conjugacy class (\cite{Nico3}). 
\end{itemize}

%

We recall that the subalgebra $B_{\varepsilon}$ 
 contains a copy of the
 coordinate ring ${\mathbb{C}}[B^-]$ of $B^-$.
 Given $b\in B^-$, let us denote
 by $m_b$ the corresponding maximal ideal  of
${\mathbb{C}}[B^-]$ and let us consider the algebra
$A_b:=B_{\varepsilon}/m_bB_{\varepsilon}$. This is a
finite-dimensional algebra with the following properties: 
\begin{theorem}\cite{3}\label{T1}
If $p$, $q\in B^-$ lie over the same element $w\in W$, then the
algebras $A_p$ and $A_q$ are isomorphic.
\end{theorem}

\begin{theorem}\cite{3}\label{T2} Let $p\in B^-\cap B\dot{w}B$ be a point
over $w\in W$ and let $A_p$ be the corresponding algebra. Assume
that $\ell$ is a good integer. Then the dimension of each irreducible
representation of $A_p$ is equal to
$\ell^{(\ell(w)+rk(1-w))/2}$. 
\end{theorem} 
\begin{corollary}\label{C1}
If $p\in B^-$ lies over $w\in  W$ and $\sigma$ is an irreducible
representation of ${\cal U}_{\varepsilon}({\mathfrak{g}})$ on a vector space
$V$ such that $\varphi(\sigma)$ is conjugated to $p$, then $\dim(V)$
is divisible by $\ell^{(\ell(w)+rk(1-w))/2}$.
\end{corollary}
\Pf See \cite[Corollary 2.9]{Nico2}.\hfill$\Box$
%

\bigskip
Theorem \ref{generale} and  Corollary \ref{C1} 
lead to the following result:
\begin{theorem}\label{main} Let $\mathfrak{g}$ be a simple
complex Lie algebra and let $\ell$ be a good integer. If $V$ is a
simple ${\cal U}_\varepsilon(\mathfrak{g})$-module 
whose associated conjugacy class $\O_V$ is spherical, then
$\ell^{\frac{1}{2} \dim \O_V}$ divides $\dim V$.\hfill$\Box$
\end{theorem} 

It was shown in \cite[\S 8]{reduction} that in order to prove the conjecture
it is enough to consider the exceptional conjugacy classes, that is the
conjugacy classes of exceptional elements.
For convenience of the reader we recall that
a semisimple element $g\in G$ is {\em exceptional} if its centralizer
in $G$ has finite centre. An element $g\in G$ is called exceptional if
its semisimple part is exceptional.
From the classification of the
semisimple exceptional elements (\cite[Lemma
7.1]{K}, \cite[\S 7]{reduction}) it follows that when $\g$ is of
classical type or of type $G_2$ all the semisimple exceptional elements are
spherical. 
The elements $\sigma_k$ and $\rho_k$ in Table 1, with $k=1,\,\ldots\left[\frac{n}{2}\right]$ 
for $\g$ of type $C_n$,  $k=2,\,\ldots\left[\frac{n}{2}\right]$ 
for $\g$ of type  $D_n$, and $k=1,\,\ldots,\,n$ for $\g$ of type $B_n$,
 are, up to central elements, representatives of all semisimple,
exceptional conjugacy classes. 
The elements appearing in Table 2
except $p_2$ and $q_3$  are, up to central elements, representatives
of all
spherical, semisimple, exceptional conjugacy classes for $\g$
exceptional type.

Using the De Concini-Kac reduction theorem (\cite[\S 8]{reduction})
we can go a bit further in
the proof of the conjecture:
\begin{corollary} Let $\g$ be of classical type or of type $G_2$ and let $s$ be a semisimple element of $G$. Then  any irreducible representation $V$ of $\Ue$ lying over
$\O_s$ has dimension divisible by $\ell^{\frac{1}{2}\dim \O_s}$.
\end{corollary} 
\proof
 An irreducible
representation of $\Ue$ lying over a semisimple element of $G$ is
either exceptional or induced by
an exceptional semisimple representation of ${\cal
U}_{\e}(\mathfrak{g}^{\prime})$ (\cite[\S 8]{reduction}). By Theorem
\ref{semisimple} the De Concini-Kac-Procesi conjecture follows for all
irreducible representations lying over semisimple elements.\hfill$\Box$

\begin{corollary}\label{larger} Let $g$ be a non-exceptional element of $G$ with Jordan decomposition
$g=su$ such that $\O_s$ and $\O_u$ are
spherical. Then any irreducible representation $V$ of $\Ue$ lying over
$\O_g$ has dimension divisible by $\ell^{\frac{1}{2}\dim \O_g}$.
\end{corollary}
\proof Since $\O_s$ is spherical $s$ can be chosen among the
non-exceptional elements in Tables 1 and 2.
The case of $G$ of type $A_n$ was dealt with in
\cite{Nico2}. Using \cite[\S 8]{reduction} we have:
\begin{enumerate}
\item $G$ is of type $C_n$ and $s=c_\lambda$.
Then $V$ is induced by an irreducible ${\cal
U}_{\e}(sp_{2n-2})$-module $V^{\prime}$ lying over the spherical, unipotent  conjugacy class of the element
 $\left(\begin{array}{c|c}U_1&U_2\\ \hline 
U_3&U_4\end{array}
\right)\in Sp_{2n-2}$;
\item $G$ is of type $C_n$ (resp. $D_n$) and $s=c$.
Then $V$ is induced by an irreducible 
 ${\cal
U}_{\e}(sl_n)$-module $V^{\prime}$ lying over the unipotent spherical conjugacy class of the element
 $A\in SL_n$ where $A$ is as in the proof of Proposition \ref{clsf}.
\item $G$ is of type $D_n$ and $s=d=\hat\tau(c)$. Then, since $C_G(c)$ is
generated by the root subgroups corresponding to the simple roots
$\alpha_1,\,\ldots,\,\alpha_{n-1}$, the centralizer of $d$ in $G$ is
generated by the root subgroups corresponding to the simple roots
$\alpha_1,\,\ldots,\,\alpha_{n-2},\,\alpha_n$ and $V$ is induced by an   
 irreducible 
 ${\cal U}_{\e}(sl_n)$-module $V^{\prime}$ lying over a spherical unipotent
conjugacy class of $SL_n$.
\item $G$ is of type $D_n$ and $s=\sigma_1$. Then  
$$u=\left(\begin{array}{cc|cc} 1&&&\\
&U_1&&U_2\\
\hline
&&1&\\
&U_3&&U_4\end{array}\right)$$ where
$u^{\prime}=\left(\begin{array}{c|c}U_1&U_2\\ \hline 
U_3&U_4\end{array}
\right)$ is a spherical unipotent element of $SO_{2n-2}$. Then
$V$ is induced by an irreducible ${\cal
U}_{\e}(so_{2n-2})$-module $V^{\prime}$ lying over $u^{\prime}$;

\item $G$ is of type $B_n$ and there are two possibilities:
\begin{enumerate}
\item $s=\rho_1$ and, as in the previous case, $V$ is induced by an irreducible ${\cal
U}_{\e}(so_{2n-1})$-module $V^{\prime}$ lying over a spherical unipotent element;
\item $s=b_\lambda$, and $V$ is induced by an irreducible ${\cal
U}_{\e}(sl_n)$-module $V^{\prime}$ lying over the unipotent spherical conjugacy class of the element
 $A\in SL_n$ where $A$ is as in the proof of Proposition \ref{clsf}.
\end{enumerate}
\item $G$ is of type $E_6$. 

In this case $s=\exp(\pi i\check{\omega}_1)$ and, since $u$ commutes with $s$, $u$ belongs to the subgroup of type $D_5$ with simple roots  $\alpha_2,\,\ldots, \alpha_6$. By \cite[\S 8]{reduction} $V$ is induced by an  irreducible 
 ${\cal
U}_{\e}(so_{10})$-module $V^{\prime}$ corresponding to the conjugacy
 class of the element $u$. Besides, the conjugacy class of $u$ in
 $D_5$ is again spherical by Remark \ref{height}. 

\item $G$ is of type $E_7$.

In this case $s=\exp(\pi i\check{\omega}_7)$ and $u$ belongs to the subgroup of type $E_6$ with simple roots  $\alpha_1,\,\ldots \alpha_6$. By \cite[\S 8]{reduction} $V$ is induced by an  irreducible 
 ${\cal
U}_{\e}(e_6)$-module $V^{\prime}$ corresponding to the unipotent
 spherical conjugacy class of the element $u$. The conjugacy class of
 $u$ in $E_6$ is again spherical by Remark \ref{height}. 
\end{enumerate}
By Theorem \ref{main} the proof is concluded.\hfill$\Box$

\begin{remark}\label{ancora-reduction}{\rm We point out that 
Corollary \ref{larger} can be generalized to a larger class of 
representations by making use of the
De Concini-Kac reduction theorem. In particular the conjecture follows
whenever the following conditions are satisfied:
\begin{enumerate}
\item $s$ lies in the identity component of
$Z(C_G(s))$; 
\item $\O_u$ is spherical.
\end{enumerate}
 
When $\g$ is of classical type condition $1.$ is equivalent to the following conditions in the corresponding
matrix groups:
\begin{itemize}
\item{\em $G=SO_{2n+1}$:} $C_G(s)^\circ$ contains no copy of type $D_k$
with $k\ge2$, i.e., if $s$ is diagonal, no submatrix of $s$ is conjugated
to $\rho_k$ with $k\ge2$; 
\item{\em $G=Sp_{2n}$:} $C_G(s)$ contains at most one copy of type
$C_k$ with $k\ge1$, i.e., if $s$ is diagonal, no submatrix of $s$ is conjugated
to $\sigma_k$ with $k\ge1$; 
\item{\em $G=SO_{2n}$:} $~C_G(s)^\circ$ contains at most one copy of type
$D_k$ with $k\ge2$, i.e., if $s$ is diagonal, no submatrix of $s$ is conjugated
to $\sigma_k$ with $k\ge2$.  
\end{itemize}

Let us notice that when $\g$ is of type $A_n$ condition $1.$ is always satisfied (\cite[Theorem
3.4]{Nico2}).}
\end{remark}

\begin{corollary} Any irreducible representation $V$ of $~{\cal U}_{\e}(sp_4)$ 
has dimension divisible by $\ell^{\frac{1}{2}\dim \O_V}$.
\end{corollary}
\proof Thanks to the De Concini-Kac reduction theorem
it is enough to consider the exceptional representations of ${\cal
U}_{\e}(sp_4)$. Since an exceptional element $g$ of $Sp_4$ is
either spherical or regular,
the
De Concini-Kac-Procesi conjecture follows from Theorem \ref{main} and
\cite[Theorem 5.1]{DC-K-P2}.\hfill$\Box$

$$$$
\end{document}